\documentclass{amsart}     

\usepackage[usenames,dvipsnames]{xcolor}                         

\usepackage{amsmath}
\usepackage{amssymb}
\usepackage{amsthm}
\usepackage{fix-cm}

\usepackage{enumerate}

\newtheorem{theorem}{Theorem}
\newtheorem{corollary}{Corollary}
\newtheorem{lemma}{Lemma}
\theoremstyle{definition}
\newtheorem{definition}{Definition}

\theoremstyle{definition}
\theoremstyle{remark}
\newtheorem{remark}{Remark}
\theoremstyle{definition}
\newtheorem{example}{Example}
\newtheorem{proposition}{Proposition}
\theoremstyle{definition}

\usepackage[utf8]{inputenc}
\usepackage{turnstile}
\usepackage{xcolor}
\usepackage{cancel}
\usepackage{ebproof}
\usepackage{float}
\usepackage[hang,flushmargin]{footmisc} 
\floatstyle{boxed} 
\restylefloat{figure}
\usepackage{pdflscape}

\usepackage{refcount}
\newcounter{fncntr}
\newcommand{\fnmark}[1]{\refstepcounter{fncntr}\label{#1}\footnotemark[\getrefnumber{#1}]}
\newcommand{\fntext}[2]{\footnotetext[\getrefnumber{#1}]{#2}}

\newcommand{\mydots}{\ifmmode\mathinner{\kern-0.1em\ldotp\kern-0.1em\ldotp\kern-0.1em\ldotp\kern-0.1em}\else.\kern-0.13em.\kern-0.13em.\fi}

\newcommand{\jmydots}{\ifmmode\mathinner{\kern-0.2em\ldotp\kern-0.2em\ldotp\kern-0.2em\ldotp\kern-0.2em}\else.\kern-0.2em.\kern-0.2em.\fi}

\newcommand{\Bdots}{\ifmmode\mathinner{\kern1mm\ldotp\kern2mm\ldotp\kern2mm\ldotp\kern1mm}\else.\kern-0.13em.\kern-0.13em.\fi}

\newcommand{\srcsize}{\@setfontsize{\srcsize}{4pt}{4pt}}

\definecolor{shadecolor}{HTML}{EEF0E7}
\definecolor{framecolor}{HTML}{892A54}
\definecolor{labelkey}{HTML}{892A54}
{\endMakeFramed}

\DeclareFontFamily{U} {MnSymbolA}{}

\DeclareFontShape{U}{MnSymbolA}{m}{n}{
	<-6> MnSymbolA5
	<6-7> MnSymbolA6
	<7-8> MnSymbolA7
	<8-9> MnSymbolA8
	<9-10> MnSymbolA9
	<10-12> MnSymbolA10
	<12-> MnSymbolA12}{}
\DeclareFontShape{U}{MnSymbolA}{b}{n}{
	<-6> MnSymbolA-Bold5
	<6-7> MnSymbolA-Bold6
	<7-8> MnSymbolA-Bold7
	<8-9> MnSymbolA-Bold8
	<9-10> MnSymbolA-Bold9
	<10-12> MnSymbolA-Bold10
	<12-> MnSymbolA-Bold12}{}

\DeclareSymbolFont{MnSyA}{U}{MnSymbolA}{m}{n}
\DeclareMathSymbol{\VDash}{\mathop}{MnSyA}{240}
\DeclareMathSymbol{\downsquig}{\mathop}{MnSyA}{163}
\DeclareMathSymbol{\rightsquig}{\mathop}{MnSyA}{160}

\title{A Cut-Free Sequent Calculus for Defeasible Erotetic Inferences}
\author[J. Millson]{Jared Millson}
\thanks{This work has benefited from the comments and suggestions of Andrzej Wi\'{s}niewski and two anynomous referees at \textit{Studia Logica}.}

\begin{document}
	
	\setcounter{page}{1}     
	
	


	
		
	\begin{center}
		\textbf{This is a pre-print of an article published in \textit{Studia Logica}. \\The final authenticated version is DOI: 10.1007/s11225-018-9839-z.}
	\end{center}
\vspace{3mm}

	 \maketitle
	
	\vspace{-.75cm}
	
	\begin{center}
		Agnes Scott College\\	
	Department of Philosophy\\
	jmillson@agnesscott.edu
	\end{center}

	\begin{abstract}
		In recent years, the effort to formalize \textit{erotetic inferences}---i.e., inferences to and from questions---has become a central concern for those working in erotetic logic. However, few have sought to formulate a proof theory for these inferences. To fill this lacuna, we construct a calculus for (classes of) sequents that are sound and complete for two species of erotetic inferences studied by Inferential Erotetic Logic (IEL): erotetic evocation and erotetic implication. While an effort has been made to axiomatize the former in a sequent system, there is currently no proof theory for the latter. Moreover, the extant axiomatization of erotetic evocation fails to capture its defeasible character and provides no rules for introducing or eliminating question-forming operators. In contrast, our calculus encodes defeasibility conditions on sequents and provides rules governing the introduction and elimination of erotetic formulas. We demonstrate that an elimination theorem holds for a version of the cut rule that applies to both declarative and erotetic formulas and that the rules for the axiomatic account of question evocation in IEL are admissible in our system.
	\end{abstract}
	


	\vspace{.5cm}
	\sloppy
	
	Most work in erotetic logic has been devoted to the relationship between questions and answers \cite{Hamblin1958,Harrah1961,Belnap1976,Karttunen1977,Hintikka1976,Groenendijk1984}. In recent years, however, the effort to formalize \textit{erotetic inferences}---i.e., inferences to and from questions---has become a central concern for those working in the subfield, e.g., \cite{Hintikka2007,Pelis2017,Wisniewski1995,Ciardelli2015a,Wisniewski2013}. Most contributions in this area deal with the semantic relations that hold among questions and statements. Although no formal study of these inferences will be complete without a treatment of the rules governing the construction of proofs, few have sought to formulate a proof theory for erotetic inferences.\fnmark{erotetic-PT} This paper aims, among other things, to fill this lacuna.
	
	\fntext{erotetic-PT}{Work on \textit{Socratic Proofs} represents a notable exception. First identified by \cite{Wisniewski2004}, Socratic Proofs are sequences of questions regarding the derivability of formulas that end with a question whose affirmative answer is, in a sense, evident. Systems for such proofs typically consist in hypersequent calculi defined over erotetically-enriched languages. Other contributions include \cite{Wisniewski2005,Leszczyska-Jasion2013}.}
	
	While several accounts of erotetic inferences have been put forward, Wi\'{s}niewski's \textsc{Inferential Erotetic Logic} (IEL) is among the most developed. There are two types of inference that IEL seeks to formalize. The first type covers cases in which one infers a question from a set of declarative statements. The second type refers to instances when one infers a question on the basis of another question and a (possibly empty) set of declaratives. According to IEL, the relation of \textit{erotetic evocation} underwrites valid instances of the first type, while that of \textit{erotetic implication} underwrites those of the second.
	
	Recently, \cite{Wisniewski2016} has attempted to axiomatize erotetic evocation in a sequent calculus defined over a classical propositional language enriched to include formal counterparts to questions---so-called \textit{erotetic formulas}. This calculus represents an important contribution to the development of a proof theory for erotetic inferences. Unfortunately, the work suffers three limitations. First, although evocative inferences are defeated when one or more answers to the evoked question are entailed by the premises (i.e., when the question is already answered by the statements that give rise to it), Wi\'{s}niewski's system does not capture the defeasible nature of this relation but instead utilizes axioms and rules that prevent defeat from arising.\fnmark{Meheus} Second, the system leaves erotetic implication completely unrepresented. Finally, there are no rules governing the introduction or elimination of erotetic formulas. The absence of such rules obscures the role that erotetic evocations play in our wider inferential economy.
	
	\fntext{Meheus}{The defeasible character of erotetic evocation \textit{has} been captured by \cite{Meheus2001}, who provides a proof theory for question evocation in the framework of adaptive logic. The present work incorporates and extends many of the insights of \cite{Meheus2001} and does so in the context of a cut-free sequent calculus.}
	
	This paper seeks to overcome each of these limitations by articulating a sequent calculus that represents species of both erotetic evocation and erotetic implication and that provides introduction and elimination rules for the interrogative operator. We follow \cite{Wisniewski2016} in limiting our target to erotetic inferences expressed in the language of classical propositional logic (CPL) extended to include erotetic formulas. Dealing with sequents formulated in an extension of CPL provides a simple and natural starting point for proof-theoretic investigations and enables us to utilize a set-based variant of Gentzen's sequent calculus for classical logic, $ \mathsf{LK} $. 
	
	Our object of inquiry is erotetic evocation and a species of erotetic implication known as \textit{regular} erotetic implication (see Definition \ref{def:reg_e-imp}). Restricting our attention to this species of erotetic implication has the benefit of yielding inferential relations that are closed under transitivity. Consequently, we are able to include a version of the cut rule that applies to both statements and questions and to show that it may be eliminated from any proof, thereby demonstrating that our system enjoys the prized property of cut-elimination exhibited by $ \mathsf{LK} $. 
	
	If the calculus proposed here is to be a genuine proof theory and not just a formal auxiliary to the semantic account of erotetic inference, it must be intelligible in terms of primitive non-semantic concepts. But unlike the natural deduction and sequent systems for CPL, which may legitimately appeal to the notion of \textit{proof} rather than that of \textit{truth-value}, it is not clear that there is an antecedently intelligible notion of \textit{proof} that applies to both statements and questions. What, then, is the appropriate interpretant for an erotetic proof theory? In answering this question, the paper does something else that is novel. It offers an informal interpretation of erotetic sequents, the proof-theoretic analogs to erotetic inferences, in terms of the normative statuses that attach to the speech acts of asking questions and making assertions. Roughly, the idea is to read sequents as saying that anyone who is entitled to assert/ask the statements/questions in the antecedent is entitled to assert/ask the statements/questions of its succedent. There are thus two primitive concepts invoked: that of being entitled to assert (a.k.a assertional entitlement) and that of being entitled to ask a question (a.k.a erotetic entitlement). This  interpretation of erotetic arguments represents an extension of the normative-pragmatic accounts of (assertional) inferences and of various non-assertional speech acts (though not of asking questions) offered by  \cite{Brandom1994} and \cite{Kukla2009}, respectively. While it is not the goal of this paper to defend the sufficiency of this interpretation, we hope that the technical results of the paper will provide fodder for such a defense.\fnmark{Diss}
	
	\fntext{Diss}{A full-blown normative-pragmatic treatment of the act of asking a question is provided by \cite{Millson2014a,Millson2014b}.}
	
	The paper proceeds as follows. In Section \ref{sec:erotetic_inferences}, we formally define the erotetic inferences that constitute the targets of our proof theory. In Section \ref{sec:SC_DR}, we introduce our sequent calculus for defeasible inferences, $ \mathsf{LK}^\mathbf{?}$. After explaining the declarative and erotetic rules of the system, we proceed to prove soundness and completeness for our target erotetic inferences---\textit{erotetic evocation} and \textit{regular erotetic implication} (Section \ref{sec:S+C}). We also prove a cut-elimination theorem for the calculus (Section \ref{subsec:cut-elim}). We then demonstrate that the axioms and primitive rules that \cite{Wisniewski2016} provides for erotetic evocation are admissible in $ \mathsf{LK}^\mathbf{?}$ (Section \ref{sec:PMC}). We conclude by exploring possible directions of future research on erotetic calculi.

	\section{Erotetic Inferences}\label{sec:erotetic_inferences}
	
	In what follows, we assume a language of CPL, $ \mathcal{L}_d $, that consists of a countable set of atoms denoted by $p,q,r,p_1 \mydots,p_n, q_1 \mydots,q_n, r_1 \mydots,r_n$, and formulas formed in the usual manner by the standard connectives $\wedge, \vee, \to, \neg$. We refer to the members of $\mathcal{L}_d $ as \textit{d-wffs}. In order to represent questions, we extend $ \mathcal{L}_d $ to $ \mathcal{L}$ by enriching the former's vocabulary with ``?,'' ``\{ ,'' ``\},'' and the comma. Well-formed expressions composed of these symbols are called \textit{e-formulas}. We let $ \mathcal{L}_e $ be the set of e-formulas, i.e., $  \mathcal{L}_e = \mathcal{L}-\mathcal{L}_d  $.
	
	We use the letters $A, B, C, D$ and $ Q, R $ as metavariables for d-wffs and e-formulas, respectively. Arbitrary formulas of either type are ranged over by $ F, G, H $. The letters $ \Gamma, \Delta, \Theta, \Sigma $ range over subsets of $ \mathcal{L} $ (i.e., sets of both declarative and erotetic formulas). Arbitrary sets of d-wffs, i.e., subsets of $ \mathcal{L}_d $ will be ranged over by $ X, \varUpsilon, Z $, while the bold letters, $ \mathbf{S}, \mathbf{T}, \mathbf{U}, \mathbf{V}$ are reserved for sets of subsets of $\mathcal{L}_d $.  Metavariables may be indexed with subscripts. All sets are assumed to be finite. We will use `$ = $' to denote both set-theoretic identity as well as syntactic equivalence of symbols.
	
	\begin{definition}[Syntax of $ \mathcal{L} $]\label{def:syntax_L}
		$ \mathcal{L} $ is the smallest set meeting the following criteria:
		\begin{itemize}
			
			\item[(i)] If $ A \in \mathcal{L}_d$, then $ A \in \mathcal{L} $.
			
			\item[(ii)] If	$ A_{1},\mydots,A_{n}  (n >1) \in \mathcal{L}_d$ and $  A_{1},\mydots,A_{n}$ are pairwise syntactically distinct, i.e., non-equiform, then $ ?\{A_{1},\mydots,A_{n}\} \in \mathcal{L}$.
			
			\item[(iii)] Nothing else is a member of $ \mathcal{L}$.
			
		\end{itemize}
	\end{definition}
	
	IEL subscribes to the \textit{set-of-answers} methodology according to which questions are identified with the set of their possible direct answers \cite{Wisniewski2013}. In keeping with this approach, we will often say that an e-formula, e.g., $?\{ A_{1},\mydots,A_{n}\} $, is \textit{based} on a particular set of d-wffs, e.g., $\{ A_{1},\mydots,A_{n}\} $, and that members of the latter set are the \textit{constituents} of the e-formula. Such comments should not mislead the reader into thinking that `?' is an operation on sets. As evidence to the contrary, note that it follows from Definition \ref{def:syntax_L} that $?\{p, \neg p \}$ and $?\{ \neg p, p\}$ are syntactically distinct e-formulas. We will often refer to e-formulas as \textit{questions} and to the d-wffs constituting them as their \textit{direct answers}. When $Q = ?\{ A_{1},\mydots,A_{n}\} $ we write $ dQ $ to denote the set of direct answers to $ Q $, i.e., $dQ =\{ A_{1},\mydots,A_{n}\} $. Since the set-notation often becomes cumbersome, we avail ourselves of some handy abbreviations.
	
	\begin{definition}[Abbreviations]\label{def:abbrev}
		We abbreviate sets and e-formulas as follows:
		
		\begin{itemize}
			\item For $ ?\{A_{1},\mydots,A_{n}\} $, we write $ ?[A_{\mid n}] $.
			\item For $ ?\{A ,\neg A\} $, we write $ ?A $.
			\item For $ \Gamma_1 \cup \mydots \cup \Gamma_n$ and $ \mathbf{S}_1 \cup \mydots \cup \mathbf{S}_n$, we write $ \Gamma_{\mid n} $ and $ \mathbf{S}_{\mid n} $, respectively.
			\item For $ \{A_{1},\mydots,A_{n}\} $, we write $ [A_{\mid n}] $.
			\item For $ \{\{A_{1}\},\mydots,\{A_{n}\}\} $, we write $ [\mathbb{A}_{\mid n}]   $.
		\end{itemize}	
	\end{definition}
	
	\subsection{Erotetic Evocation}\label{subsec:e-evocation+implication}
	
	The concepts of erotetic evocation and implication are explicated by means of standard, single-conclusion entailment or \textit{sc-entailment} and multiple-conclusion entailment or \textit{mc-entailment}. IEL defines both quite broadly so as to accommodate a diversity of underlying languages. Since our focus here is only on instances of these relations that are formulated in CPL, we will simply use entailment in that logic for the first semantic relation and will define mc-entailment accordingly.  Thus, we assume a Boolean valuation, $ v $ that assigns a truth value, \textbf{1} or \textbf{0}, to each atom in $ \mathcal{L}_d $ and extends to all d-wffs in the usual manner by using the Boolean functions corresponding to the connectives. By $ v(A) = \mathbf{1} $ we mean ``d-wff A is true under valuation $ v $.'' The definitions for sc- and mc-entailment are straightforward.

	\begin{definition}[SC-Entailment, $ X \vDash B $]\label{def:sc-entailment}
		$ X \vDash B $ if and only if there is no valuation $ v $ such that $ v(A) = \mathbf{1}$ for every $ A \in X $ and $ v(B) = \mathbf{0}$.
	\end{definition}
	
	\begin{definition}[MC-Entailment, $ X \VDash \varUpsilon $]\label{def:mc-entailment}
		$ X \VDash \varUpsilon $ if and only if there is no valuation $ v $ such that $ v(A) = \mathbf{1}$ for every $ A \in X $ and $ v(B) = \mathbf{0}$ for every $ B \in \varUpsilon $.
		
	\end{definition}
	
	Before defining our target erotetic inferences, we introduce two useful concepts. The first is the property of \textit{soundness} that a question possesses when at least one of its direct answers is true. In CPL, the soundness of a question may defined as follows.
	
	\begin{definition}[Soundness in CPL]
		A question, $ Q $ is \textit{sound under v} iff $v(\bigvee dQ) = \mathbf{1}$.
	\end{definition}
	
	The second concept captures the relation between a set of declaratives $ X $ and a question  $ Q $ such that it is impossible that all the d-wffs in $ X $ are true, but no direct answer to $ Q $ is true. When this relation holds, we say that $ Q $ is \textit{sound relative} to $ X $.
	
	\begin{definition}[Relative Soundness]
		A question $ Q $ is sound relative to $ X $ iff $ X\VDash dQ $.
	\end{definition}

	\subsection{Erotetic Implication}\label{subsec:e-evocation+implication2}
	We now define erotetic evocation.
	
	\begin{definition}[Erotetic Evocation (EE), $ \mathbf{E}(X, Q) $]\label{def:erotetic_evocation}
		A set of d-wffs, $ X $ evokes a question, $ Q $, i.e., $ \mathbf{E}(X, Q)  $ if and only if 
		\begin{itemize}
			\item[(i)] $ X \VDash dQ $, and
			\item[(ii)] $\neg\exists A \in dQ (X \vDash A) $.
		\end{itemize}
	\end{definition}
	
	The first condition for erotetic evocation requires that an evoked question, $ Q $ has a true direct answer if the evoking set consists of truths, in other words, that $ Q $ is \textit{sound relative} to $ X $. The second condition ensures that evoked questions are informative relative to $ X $, i.e., that the evoking premises do not entail an answer. The latter makes EE a defeasible relation in the sense that a question that is evoked by $ X $ will not be evoked by a proper superset of $ X $ that includes formulas that sc-entail a direct answer to $ Q $.
	
	We now define erotetic implication and the species of that relation that constitutes our second target. In keeping with the literature on IEL, this definition uses `$ \subset $' to denote the \textit{proper} inclusion of sets.
	
	\begin{definition}[Erotetic Implication (EI), $ \mathbf{Im}(Q, X, Q_1) $]\label{def:e-imp}
		A question $ Q $ implies a question $ Q_1 $ on the basis of a set of d-wffs $ X $, i.e., $ \mathbf{Im}(Q, X, Q_1) $, if and only if
		\begin{itemize}
			\item[(i)] $ \forall A \in dQ (X \cup \{A\} \VDash dQ_1 )$, and
			\item[(ii)] $ \forall B \in dQ_1\:\,\exists \varUpsilon\subset dQ  \text{\:s.t.\:} \varUpsilon \neq \emptyset \:(X \cup \{B\} \VDash \varUpsilon) $.
		\end{itemize}
	\end{definition}
	
	The first condition for erotetic implication says that the implying question, together with a set of auxiliary d-wffs, $ X $, must mc-entail the set of direct answers to the implied question. This condition ensures that the truth/soundness of the declarative/erotetic premises is \textit{transmitted} to the erotetic conclusion. In other words, if all of  $ X$ are true, $ Q $ is \textit{sound}, and $ \mathbf{Im}(Q, X, Q_1) $, then $ Q_1 $ is \textit{sound}. 
	
	The second condition guarantees that $ Q_1 $ is \textit{cognitively useful} relative to $ X $ and $ Q $. By `cognitively useful,' we mean that each direct answer to the implied question should, together with the declarative premises, narrow down the class of possible direct answers to the implying question.
	
	We will focus on a special kind of erotetic implication, namely, \textit{regular} erotetic implications. Such a restriction is desirable because regular erotetic implication is transitive. This means that our sequent calculus will enjoy a cut rule that applies to e-formulas as well as to d-wffs.
	
	\begin{definition}[Regular Erotetic Implication (R-EI), $ \mathbf{Im}_{r}(Q, X, Q_1) $]\label{def:reg_e-imp}
		A question $ Q $ regularly implies a question $ Q_1 $ on the basis of a set of d-wffs $ X $, i.e., $ \mathbf{Im}_{r}(Q, X, Q_1) $, if and only if
		
		\begin{itemize}
			\item[(i)] $ \forall A \in dQ (X \cup \{A\} \VDash dQ_1 )$, and
			\item[(ii)] $ \forall B \in dQ_1\:\,\exists A \in dQ (X \cup \{B\} \vDash A) $.
		\end{itemize}
	\end{definition}
	
	\begin{corollary}[Transitivity of R-EI]\label{coro:transivity_R-EI}
		If $ Q $ regularly implies $ Q_1 $ on the basis of $ X $, and $ Q_1 $ regularly implies $ Q_2$ on the basis of $ X $, then $ Q $ regularly implies $ Q_2$ on the basis of $ X $.
	\end{corollary}
	
	When the set of auxiliary d-wffs is empty, the relation of erotetic implication is said to be \textit{pure}.
	
	\begin{definition}[Regular Pure Erotetic Implication (RP-EI), $ \mathbf{Im}_{rp}(Q, Q_1) $]\label{def:r-p_e-imp}
		A question $ Q $ regularly and purely implies a question $ Q_1 $, i.e., $ \mathbf{Im}_{rp}(Q, Q_1) $, if and only if
		
		\begin{itemize}
			\item[(i)] $ \forall A \in dQ (\{A\} \VDash dQ_1 )$, and
			\item[(ii)] $ \forall B \in dQ_1\:\,\exists A \in dQ (\{B\} \vDash A) $.
		\end{itemize}
	\end{definition}
	
	As we shall demonstrate, proper subsets of provable sequents in our calculus are sound and complete with respect to the definitions of R-EI and RP-EI.
	
	\section{A Sequent Calculus for Defeasible Reasoning}\label{sec:SC_DR}
	
	We will now present our sequent calculus for defeasible reasoning. We employ the standard notation for sequents composed of sets of formulas---e.g., $\Gamma, A \vdash B, \Delta $. Formulas on the left side of the turnstile are called the \textit{antecedent}; on the right side they are called the \textit{succedent}. Commas in the antecedent are read conjunctively and those on the right are read disjunctively. The formula with the connective in a rule is the \textit{principal} formula of that rule, and its components in the premises are the \textit{active} formulas. The remaining elements of sequents are referred to as \textit{side formulas}.
	
	Sequents in our calculus represent rules of inference. As such rules are inherently normative, we interpret sequents as representing a relationship that preserves a particular normative status, namely, entitlement. Thus, a provable sequent tells us that anyone who is entitled to the formulas in the antecedent is entitled to at least one of those in the succedent. Of course, `entitlement to a formula' is at best elliptical for `entitlement to \textit{do something} with that formula'. So, our reading of sequents requires a further, \textit{pragmatic} layer of interpretation. In particular, a provable sequent says that anyone entitled to \textit{assert} the declaratives and \textit{ask} the questions in the antecedent is entitled to \textit{assert} or \textit{ask} at least one of the declaratives/questions in the succedent. This is the normative-pragmatic interpretation that underwrites our proof theory. While sequents represent rules of entitlement-preserving inference, the rules of our calculus are meta-rules of inference---i.e., they tell us what rules of inference are legitimate given the axioms.

	The sequents in our calculus depart from the standard form in two respects: just below our turnstile we add a set of sets of d-wffs, $\mathbf{S}$, called a \textit{defeater set} and to the far left of the turnstile we add a set of formulas, $\Sigma $, called a \textit{background set}.$$ \Sigma\mid\Gamma\sststile{\mathbf{S}}{}\Delta $$
	
	Defeater sets `describe' situations in which one is not permitted to draw the inference represented by the sequent. In line with our normative-pragmatic reading of sequents, we can think of the elements of (members of) the defeater set as those declaratives entitlement to which prohibits entitlement to draw the inference represented by the sequent.
	
	We interpret background sets as consisting of information to which a reasoner is beholden when she draws an inference, but which does not serve as a premise and from which the conclusion is not said to follow. Having such a device in our formalism enables us to capture an important aspect of defeasibility, namely, that the introduction of new information may jeopardize prior inferential commitments even when that information does not serve as fodder for new inferences. These sets also serve a technical role in our system since they often expand in the course of a derivation, keeping a \textit{record}, so to speak, of those formulas that `disappear' from antecedents of premises. (See $\vee\vdash, \to\vdash, \vdash\to, $ and $\neg\vdash$ in Figure 1).
	As we shall see, the latter feature proves crucial to the implementation of a Gentzen-style normalization procedure and to the establishment of a cut-elimination theorem.

	\begin{definition}[Defeater sets, Background sets, Defeasible Sequents]\label{def:def_sequents}
		Defeater sets are finite, possibly empty sets of d-wffs representing statements, entitlement to which is incompatible with an entitlement to draw the inference represented by the sequent. Background sets are finite, possibly empty sets of formulas that represent the informational context of the inference. A defeasible sequent is a standard sequent in which a \textit{background set} ($ \Sigma $) occurs on the left-hand side of the antecedent and a \textit{defeater set} ($ \mathbf{S} $) occurs just below the turnstile:\fnmark{Piazza-Control_Sets}
		$\Sigma\mid\Gamma\sststile{\mathbf{S}}{}\Delta $. When no background sets have been specified (i.e., $\Sigma = \emptyset$) we write: $ \cdot\mid\Gamma\sststile{\mathbf{S}}{}\Delta$. When the antecedent is empty we write: $ \Sigma\mid\cdot\sststile{\mathbf{S}}{}\Delta$.
		
	\end{definition}
	
	\fntext{Piazza-Control_Sets}{ In general, the respective roles played by  \textit{control sets} and \textit{soundness} in Piazza and Pulcini's $\mathsf{LK}^\mathcal{S} $ are played by \textit{defeater sets} and \textit{undefeatedness} in our system. Aside from various technical differences, the crucial distinction between the two approaches is that for  \cite{Piazza2017} the occurrence of a disjunctive formula in the antecedent renders the sequent defeated (resp. unsound) if either of the disjuncts occurs in the sequent's defeater set (resp. control set), while in our system, the sequent is only  defeated if \textit{both} disjuncts are present in the defeater set. We believe that the latter property captures a more intuitive, less cautious conception of defeat. In order to realize this conception in our definitions, we found it necessary to deploy substantially different operations and have re-named the resultant concepts so as to avoid confusion.  With this said, Definition \ref{def:proof} and Lemma \ref{lem:defeatrel} are taken over from \cite{Piazza2017} with little modification.}
	
	\begin{definition}[Derivation, Proof, Paraproof]\label{def:proof}
		If a defeasible sequent, $ \Sigma\mid\Gamma\sststile{\mathbf{S}}{}\Delta $, occurs at the root of a finitely branching tree, $ \pi $, whose nodes are defeasible sequents recursively built up from axioms by means of the rules of $ \mathsf{LK}^\mathbf{?} $, then $ \pi $ is said to be a \textit{derivation} of $ \Sigma\mid\Gamma\sststile{\mathbf{S}}{}\Delta $. If, additionally, each sequent in $ \pi $ is undefeated, then $ \pi $ is said to be a proof of $ \Sigma\mid\Gamma\sststile{\mathbf{S}}{}\Delta $ in $ \mathsf{LK}^\mathbf{?} $,\, otherwise $ \pi $ is called a \textit{paraproof} of $ \Sigma\mid\Gamma\sststile{\mathbf{S}}{}\Delta $.
	\end{definition}
	
	The rules for $ \mathsf{LK}^\mathbf{?}$ in Figure 1 
	are designed to generate trees that preserve derivability \textit{downward} and undefeatedness \textit{upward}. When read bottom-up, the rules say that, ``If [the conclusion sequent] is undefeated, then so is/are [the premise sequent/s].'' Alternatively, the rules may be read top-down as permitting the conclusion, given the premises, so long as the conclusion is not defeated. On either reading, the rules are formulated to ensure that a cut-free derivation whose end-sequent is undefeated will contain only undefeated sequents throughout. This property is important both for the establishment of cut-elimination and, more primitively, for the reason that proofs in a system of defeasible inference should not contain defeated sequents. 
	
	With the exception of $\mathsf{BE}, \mathsf{DE}, ?\vdash_1, ?\vdash_2, \vdash ?_1,$ and $ \vdash ?_2,  $ rules that apply to declarative sequents are based on \cite{Ketonen1944}'s alternative to Gentzen's $ \mathsf{LK} $ in which rules permit independent contexts (i.e. are multiplicative) and sequents are finite constructions of sets of formulas, rather than constructions of sequences or multisets. Set-based calculi of this sort obviate the structural rules for contraction and exchange and offer a simpler proof system. As we shall see, the structural similarity between the d-wff rules in our system and those found in set-based versions of $ \mathsf{LK} $ is of great help in establishing our theorems for soundness and completeness.
	
	We now define the all-important property of defeat.
	
	\begin{definition}[$ \mathcal{E}(\cdot)$]\label{def:EClo}
		Let $ \mathcal{E}(\Gamma) $ denote the set that results from removing all e-formulas from $ \Gamma $ and replacing them with a disjunction of their declarative constituents, i.e.
		$$  \mathcal{E}(\Gamma) = \{A: A\in\Gamma\}\cup\{A_1\vee\ldots\vee A_n: \:?\{A_1,\ldots,A_n\} \in \Gamma\}$$
		
		\noindent Since the $ \mathcal{E}(\cdot) $ operation replaces questions with disjunctions of their direct answers, we say that $ \mathcal{E}(\Gamma) $ is the \textit{declarativization} of $ \Gamma $.
	\end{definition}
	
	\begin{definition}[Defeat]\label{def:defeat} A sequent is defeated just in case it is possible to derive a sequent in $ \mathsf{LK^{?}} $  that has empty background and defeater sets, an antecedent that combines the decarativized background and antecedent sets of the original sequent, and a succedent that is a member of the original sequent's defeater set. We define \textit{defeat} formally as follows. 
		
		Let $\Sigma\mid\Gamma\sststile{\mathbf{S}}{\mathcal{D}_{\mathsf{LK^{?}}}}\Delta$ stand for the fact that $\Sigma\mid\Gamma\sststile{\mathbf{S}}{}\Delta$ is derivable in $ \mathsf{LK^{?}} $.
		$$   \Sigma\mid\Gamma\sststile{\mathbf{S}}{}\Delta\text{ \: is defeated \textit{iff} \:} \exists X\in\mathbf{S} \text{\, such that \,} \cdot\mid\mathcal{E}(\Sigma\cup\Gamma)\sststile{\emptyset}{\mathcal{D}_{\mathsf{LK^{?}}}} X. $$
	\end{definition}
	
	The results in section \ref{subsec:LK-S?d} show that this is equivalent to saying that a sequent in $ \mathsf{LK^{?}} $ is defeated just in case the declarativized antecedent-\textit{cum}-background set mc-entails a member of the defeater set. But while this fact may aid the reader in identifying defeated sequents, it does not constitute its formal definition. Rather, we define defeat in terms of derivability in $ \mathsf{LK^{?}} $ and provability in terms of derivability and defeat---thereby ensuring that all the resources needed for proof search are internal to the system itself.

	\begin{example}\label{ex:con_in_ant}
		$ \quad\cdot\mid p\wedge q \sststile{ \{\{p \}\} }{} \Delta $ is defeated.
	\end{example}
	
	\begin{example}\label{ex:con_in_def}
		$ \quad\cdot\mid p\sststile{\{\{p \wedge q\}\}}{} \Delta $ is undefeated.
	\end{example}
	
	\begin{example}\label{ex:dis_in_def}
		$ \quad\cdot\mid p \sststile{\{\{p\vee q \}\}}{} \Delta $ is defeated.
	\end{example}
	
	\begin{example}\label{ex:dis_in_ant}
		$ \quad\cdot\mid p \vee q \sststile{\{\{p\}\}}{} \Delta $ is undefeated.
	\end{example}
	
	\begin{example}\label{ex:Q_in_ant}
		$ \quad\cdot\mid\: ?\{p, \neg p\} \sststile{\{\{p\}\{\neg p\}\}}{} \Delta $ is undefeated.
	\end{example}
	
	\begin{example}\label{ex:Q_in_ant-def}
		$ \quad\cdot\mid\: ?\{p, \neg p\} \sststile{\{\{p, \neg p\}\}}{} \Delta $ is defeated.
	\end{example}
	
	The definition of \textit{defeat} is designed to redeem certain intuitions about defeasible inferences. For instance, if we know that $ A $ defeats an inference, then we ought to reject that inference if we are entitled to $ A\wedge B $ (Example \ref{ex:con_in_ant}).  Conversely, if we know that $ A \wedge B $ defeats an inference, but are ignorant as to whether either conjunct by itself defeats it, then we should not abandon the inference if we are merely entitled to one of the conjuncts (Example \ref{ex:con_in_def}). Similarly, if $ A\vee B $ defeats an inference, then the inference is defeated if we are entitled to either of the disjuncts (Example \ref{ex:dis_in_def}). Finally, if $ A $ defeats an inference, then entitlement to $ A\vee B $ need not force us to abandon that inference, since entitlement to $ B $ would preserve the propriety of the inference (Example \ref{ex:dis_in_ant}). 
	
	When it comes to inferences with questions among their premises, defeat should occur when the premises entail at least one of their direct answers, for there is a strong sense in which being entitled to ask a question (or better, to conduct an inquiry) is incompatible with possession of its answer (resolution). Furthermore, the fact that defeater sets are sets of sets of d-wffs enables us to distinguish between situations in which entitlement to at least one particular direct answer prohibits entitlement to ask the question (Example \ref{ex:Q_in_ant}) and those in which entitlement to any answer is foreclosed, i.e., when the question is unanswerable in the context (Example \ref{ex:Q_in_ant-def}). As we shall see, by exploiting the dynamics of the defeat mechanism, the calculus is able to represent the defeasible character of erotetic evocation.
	
	We note that sequents with inconsistent antecedent-\textit{cum}-background sets are defeated whenever their defeater sets are nonempty.
	
	\begin{lemma}\label{lem:inconsistent}
		If $ \mathcal{E}(\Sigma\cup\Gamma) $ is inconsistent and $ \mathbf{S} \neq\emptyset $, then $ \Sigma\mid\Gamma\sststile{\mathbf{S}}{}\Delta $ is defeated. 
		\begin{proof}
			According to Definition \ref{def:EClo} and the rules for $ \mathsf{LK^{?}} $, if $ \mathcal{E}(\Sigma\cup\Gamma) $ is inconsistent, then $\cdot\mid\mathcal{E}(\Sigma\cup\Gamma)\sststile{\emptyset}{\mathcal{D}_{\mathsf{LK^{?}}}} X$ for any $ X$. By hypothesis, $\exists X\in\mathbf{S} $. Thus, according to Definition \ref{def:defeat}, the sequent is defeated no matter what the defeater set happens to contain.
		\end{proof}	
	\end{lemma}
	
	\begin{definition}[Compatibility, $ \succsim $]\label{def:compat}
		When a set of formulas, $ \Gamma $, and a set of sets of d-wffs, $ \mathbf{S} $, fail to meet the conditions of defeat, they are said to be \textit{compatible}. We use the symbol `$ \succsim $' to denote this relationship, which we define formally as follows:
		$$ \Gamma\succsim\mathbf{S} \:\textit{iff}\:\,\neg\exists X\in\mathbf{S}\text{\, such that \,} \cdot\mid\mathcal{E}(\Sigma\cup\Gamma)\sststile{\emptyset}{\mathcal{D}_{\mathsf{LK^{?}}}} X.$$
		When this relation fails to hold, we say that $ \Gamma $ and $ \mathbf{S} $ are \textit{incompatible} and express this by writing: $ \Gamma\not\succsim\mathbf{S} $.
	\end{definition}
	
	We now look at a lemma that aids understanding of the rules in our calculus.
	
	\begin{lemma}\label{lem:defeatrel}
		\begin{enumerate}
			\item If $\Gamma\cup\Delta \succsim \mathbf{S}$ and $ \mathbf{T} \subset \mathbf{S} $, then $ \Delta \succsim \mathbf{T}$.
			\item $ \Gamma \cup \{A \wedge B\} \succsim \mathbf{S}$ \,\, iff\,\,  $ \Gamma \cup \{A, B\} \succsim \mathbf{S}$
			\item $ \Gamma \cup \{A \vee B\} \succsim \mathbf{S}$ \,\, iff\,\,  $ \Gamma \cup \{A\} \succsim \mathbf{S}$ \,or\, $ \Gamma \cup \{B\} \succsim \mathbf{S}$
			\item $ \Gamma \cup \{\neg\neg A\} \succsim \mathbf{S}$ \,\,iff\,\,  $\Gamma \cup \{A\} \succsim \mathbf{S}$
			\item $ \Gamma\,\cup \,?\{A_1,\mydots,A_n\} \succsim \mathbf{S}$ \,\, iff\,\,  $ \Gamma \cup \{A_1\} \succsim \mathbf{S}$ \,or\ldots or\, $ \Gamma \cup \{A_n\} \succsim \mathbf{S}$
		\end{enumerate}
		
		\begin{proof}
			For the first sub-lemma, suppose for \textit{reductio} that $ \Delta\not\succsim \mathbf{T}$. From Definition \ref{def:compat} it follows that $ \cdot\mid\mathcal{E}(\Gamma\cup\Delta)\sststile{\emptyset}{\mathcal{D}_{\mathsf{LK^{?}}}} X$  for some $ X \in \mathbf{T} $ and since $  \mathbf{T} \subset \mathbf{S} $, it follows that $ \Gamma\cup\Delta\not\succsim\mathbf{S}$, contradicting the hypothesis.  Each of the remaining sub-lemmas follows directly from Definitions \ref{def:defeat} and \ref{def:compat}.
		\end{proof}
	\end{lemma}
	
	\subsection{The declarative rules for $ \mathsf{LK}^\mathbf{?}$}\label{subsec:LK-S?d}
	
	For clarity, we will refer to the axioms and logical rules for d-wffs of $ \mathsf{LK}^\mathbf{?} $, as well as to the structural and cut rules whose active and principle formulas are exclusively d-wffs, as the \textit{declarative rules}, \textit{d-wff rules} or simply $ \mathsf{LK_d}^\mathbf{?} $. As we mentioned above, $ \mathsf{LK_d}^\mathbf{?} $ is based on Ketonen's formulation of $ \mathsf{LK} $---a formulation whose primary characteristics are multiplicative rules for connectives and set-based sequents. In fact, if we ignore the background and defeater set notation and as well as the rules that govern their manipulation, i.e., we restrict ourselves to $\mathsf{LK_d}^?\backslash\{\mathsf{BE} \cup \mathsf{DE}\} $, then we have a system that is structurally similar to $ \mathsf{LK^{\{\}}} $ in \cite[p. 64]{Bimbo2014} and $ \mathbf{Gcl^*} $ in \cite[p. 7]{Poggiolesi2011}, and that is isomorphic to $ \mathsf{LK_0} $ in \cite[p. 711]{Piazza2016}.\fnmark{LK0} This means that we can avail ourselves of the normalization procedure for set-based formulations of $ \mathsf{LK} $ and for $ \mathsf{LK_0} $ in particular.
	
	\fntext{LK0}{Aside from construing sequents as constructions of sets, and hence rendering rules for contraction and exchange redundant, $ \mathsf{LK_0} $ is distinguished by its use of atomic singletons for axioms, i.e., $ p\vdash p $.}
	
	Of course, the rules $ \mathsf{DE} $ and $ \mathsf{BE} $ are not superfluous. They play an important role in cut-elimination, as shown in \cite{Piazza2017}. They also help the system capture important episodes in defeasible reasoning---namely, instances when one's knowledge of potential defeaters ($ \mathsf{DE} $) or relevant features of the context ($ \mathsf{BE} $) expands.
	
	The rules $ \vee\vdash, \to\vdash, \vdash\to, $ and $ \vdash\neg $, preserve undefeatedness upwards by exploiting the background set device. To illustrate this property, we focus on $ \vdash\neg $. This rule adds the (active) antecedent of the premise to the background set of the conclusion. Such behavior is of a piece with the explanation given for background sets above---they act as a kind of \textit{record} of those formulas that appear on the left-hand side of a turnstile in the premises but not in the conclusion.\fnmark{backsetcut} An informal interpretation of the rule (read upwards) can be given as follows: if one is entitled to (assert) $ \neg A $ while $ A $ is in one's background set of entitlements, then one is entitled to whatever follows from $ A $ once it has been removed from that background set and entitlement to its negation has been renounced. The presence of $ A $ in the background set of the conclusion thus ensures that the premises are undefeated whenever the conclusion is. 
	
	\fntext{backsetcut}{By this same reasoning, the $ cut $ rule also ought to add the active formula in its premises to the background set of the conclusion. However, the $ cut  $ rule is exempted on the grounds that such a rule is just a statement about the conditions under which information may be \textit{removed} from a proof.}
	
	\begin{lemma}\label{lem:upundefeat}
		Any cut-free derivation in $ \mathsf{LK_d}^\mathbf{?} $ is a proof if and only if its end-sequent is undefeated, i.e., undefeatedness is preserved upwards in cut-free proofs in $ \mathsf{LK_d}^\mathbf{?} $.
		
		\begin{proof}
			For all of the rules in $ \mathsf{LK_d}^\mathbf{?} $, the undefeatedness of the premises follows directly from that of the conclusion by way of Lemma \ref{lem:defeatrel}.
		\end{proof}
	\end{lemma}

	\begin{lemma}\label{lem:cutelim}
		Any sequent that is provable in $ \mathsf{LK_d}^\mathbf{?} $ has a cut-free proof.
		
		\begin{proof}
			The normalization procedure for   $ \mathsf{LK_d}^\mathbf{?} $ is just that of $ \mathsf{LK_0} $, \textit{modulo} the rules $ \mathsf{BE} $ and $ \mathsf{DE} $. But since the latter do not generate occasions for \textit{cut}, their addition is innocuous. The main challenge posed by proofs of defeasible sequents is the fact that as \textit{cut} is permuted upward in reducing derivations, the accumulation of antecedents and background sets occasionally generates defeated sequents in the derivation that were not present in the original one. Here is an example:
			$$
			\begin{prooftree}[separation = .5em,template = \small$\inserttext$ ]
			\Hypo{\pi_1}
			\Ellipsis{}{ p \mid\Gamma\sststile{\{\{r\}\}}{}q}
			\Hypo{\pi_2}
			\Ellipsis{}{\cdot\mid\Gamma'\sststile{\mathbf{S}}{} r}
			\Infer2[$ \scriptstyle{\vdash\wedge} $]{ p\mid\Gamma', \Gamma\sststile{\{\{r\}\}\,\cup\,\mathbf{S}}{}q\wedge r}
			\Hypo{\pi_3}
			\Ellipsis{}{\cdot \mid q, r\sststile{\mathbf{T}}{}\Delta}
			\Infer1[$ \wedge\vdash $]{ \cdot\mid q\wedge r\sststile{\mathbf{T}}{} \Delta}
			\Infer2[$ \scriptstyle{cut} $]{ p\mid \Gamma',\Gamma\sststile{\{\{r\}\}\,\cup\,\mathbf{S}\,\cup\,\mathbf{T}}{} \Delta}
			\end{prooftree}
			$$
			$$ \downsquig $$
			$$
			\begin{prooftree}[separation = 1em,template = \small$\inserttext$ ]
			\Hypo{\pi_1}
			\Ellipsis{}{ p \mid\Gamma\sststile{\{\{r\}\}}{}q}
			\Hypo{\pi_3}
			\Ellipsis{}{\cdot \mid q, r\sststile{\mathbf{T}}{}\Delta}
			\Infer2[$ \scriptstyle{cut} $]{ p\mid\Gamma, r\sststile{\{\{r\}\}\,\cup\,\mathbf{T}}{}\Delta}
			\Hypo{\pi_2}
			\Ellipsis{}{\cdot\mid\Gamma'\sststile{\mathbf{S}}{} r}
			\Infer2[$ \scriptstyle{cut} $]{ p\mid \Gamma',\Gamma\sststile{\{\{r\}\}\,\cup\,\mathbf{S}\,\cup\,\mathbf{T}}{} \Delta}
			\end{prooftree}
			$$
			Notice that the defeated sequent $ p\mid\Gamma, r\sststile{\{\{r\}\}\,\cup\,\mathbf{T}}{}\Delta $ does not appear in the first proof but does in the second. Thus, we have moved from a proof to a paraproof. To be certain that the elimination of \textit{cut} from the final reduction does not produce a paraproof, we must know that the cut-free derivation of an undefeated end-sequent is a proof, and we \textit{do} know this in virtue of Lemma \ref{lem:upundefeat}. Therefore, the normalization procedure for  $ \mathsf{LK} $ may be applied to proofs in  $ \mathsf{LK_d}^\mathbf{?} $.
		\end{proof}	
	\end{lemma}
	
	With the establishment of cut-elimination for the d-wff proofs of $ \mathsf{LK}^\mathbf{?}$, we can see that the structural isomorphism between $ \mathsf{LK_d}^\mathbf{?}\backslash\{\mathsf{BE} \cup \mathsf{DE}\}  $ and $ \mathsf{LK_0} $ (minus background and defeater set notation) gives $ \mathsf{LK^?} $ some important properties. 

	\begin{definition}\label{def:vdashlk}
		Let $ X\vdash_\mathsf{LK_0}\Upsilon $ stand for the fact that $ X\vdash\Upsilon $ is provable in $ \mathsf{LK_0} $.
	\end{definition}
	
	\begin{lemma}\label{lem:soundproof-LK}
		If $ \Sigma\mid X\sststile{\mathbf{S}}{} \varUpsilon$ is provable in $ \mathsf{LK_d}^\mathbf{?}$, then $X \vdash_\mathsf{LK_0} \varUpsilon$.
		\begin{proof}
			From Definition \ref{def:proof} and the isomorphism between the rules of $ \mathsf{LK_d}^?$ and those of $ \mathsf{LK_0} $.
		\end{proof}
	\end{lemma}
	
	\begin{lemma}\label{lem:proof-LK}
		There is a proof of $ \Sigma\mid X\sststile{\emptyset}{} \varUpsilon$ in $ \mathsf{LK_d}^\mathbf{?}$ iff $ X \vdash_\mathsf{LK_0} \varUpsilon$.
		\begin{proof}
			From Definition \ref{def:defeat} it follows that a sequent with an empty defeater set is provable just in case it is derivable in $ \mathsf{LK^{?}} $. Call  $ X\vdash\Upsilon $ the $ \mathsf{LK} $\textit{-analogue} of $ \Sigma\mid X \sststile{\mathbf{S}}{}\Upsilon $. Since a sequent is derivable in $ \mathsf{LK^{?}} $ just in case its $ \mathsf{LK} $\textit{-analogue} is provable in $\mathsf{LK_0} $, we obtain our result from Lemma \ref{lem:soundproof-LK}.
		\end{proof}
	\end{lemma}
	
	Lemma \ref{lem:proof-LK} implies that the consequence relation encoded by $\sststile{\emptyset}{}  $ in $ \mathsf{LK_d}^\mathbf{?} $ corresponds to the consequence relation encoded by $ \mathsf{LK_0} $. The former is of course not the only consequence relation encoded by $ \mathsf{LK_d}^\mathbf{?} $, let alone by $ \mathsf{LK}^\mathbf{?}$. Since consequence relations in the latter are individuated by defeater sets, there are countably infinitely many potential consequence relations in $ \mathsf{LK}^\mathbf{?}$. Naturally, many of these will not be well-behaved by standard criteria. Nonetheless, \textit{classes} of these consequence relations will have important properties, and, as we shall show, some of these are sound and complete with respect to our target erotetic relations. Before proceeding to demonstrate this, we draw out some of the semantic implications of the above lemmas.
	
	Since $ \mathsf{LK_0} $ is sound and complete for CPL, it follows that if $X\vdash_\mathsf{LK_0} \Upsilon$, then $ X $ mc-entails $ \Upsilon $, i.e.
	
	\begin{lemma}[$ \vdash_\mathsf{LK_0} \Leftrightarrow \VDash$]\label{lem:LK_sound+complete}
		$ X \vdash_\mathsf{LK_0} \varUpsilon  $ iff $ X \VDash \varUpsilon $
	\end{lemma}
	
	From Lemmas \ref{lem:soundproof-LK} and \ref{lem:LK_sound+complete}, we can derive:
	
	\begin{lemma}\label{lem:LKSd_sound}
		If $ \Sigma\mid X\sststile{\mathbf{S}}{} \varUpsilon$ is provable in $\mathsf{LK_d}^\mathbf{?} $  then $ X \VDash \varUpsilon $, for arbitrary $ \Sigma $ and $ \mathbf{S} $.
	\end{lemma}
	
	From Lemmas \ref{lem:proof-LK} and \ref{lem:LK_sound+complete}, we can derive:
	
	\begin{lemma}\label{lem:LKSd_sound-complete}
		$ \Sigma\mid X\sststile{\emptyset}{} \varUpsilon$ is provable in $\mathsf{LK_d}^\mathbf{?} $ iff $ X \VDash \varUpsilon  $ for an arbitrary $ \Sigma $.
	\end{lemma}
	
	
	Recall that proofs in $\mathsf{LK}^\mathbf{?} $ must contain undefeated sequents. It follows that, so long as a defeasible sequent contains a nonempty defeater set, it is possible that there will be no proof of it in $\mathsf{LK_d}^\mathbf{?} $, despite the fact that the antecedent mc-entails the succedent. As a system for defeasible reasoning, it is to be expected that certain patterns of otherwise classically valid reasoning simply will not be permitted in $\mathsf{LK}^\mathbf{?} $. Thus, Lemma \ref{lem:LKSd_sound} only says that $\mathsf{LK_d}^\mathbf{?} $ is sound with respect to mc-entailment. However, if the defeat mechanism is bypassed in virtue of the fact that the defeater set is empty (Lemma \ref{lem:LKSd_sound-complete}), the resulting class of sequents will be both sound and complete.
	
	Finally, these insights allow us to redeem our earlier claim that a sequent is defeated just in case the declarativized antecedent-\textit{cum}-background set mc-entails a member of the defeater set.
	
	\begin{corollary}\label{coro:defeatasclassical}
		$ \Sigma \mid \Gamma \sststile{\mathbf{S}}{} \Delta $ is defeated iff $\exists X\in\mathbf{S}$  such that $ \mathcal{E}(\Sigma\cup\Gamma) \VDash X $.
	\end{corollary}
	
	\subsection{The erotetic rules for $ \mathsf{LK}^\mathbf{?}$}\label{subsec:LK-S?e}
	
	We now turn to the right and left rules for e-formulas in  $ \mathsf{LK}^\mathbf{?}$ (i.e., $ \vdash?_1 $,$ \vdash?_2 $, $ ?\vdash_1 $, and $ ?\vdash_2$). Since right and left rules for connectives in a sequent calculus correspond, respectively, to the introduction and elimination rules in natural deduction, we will say that the rules  $ \vdash?_1 $ and $ \vdash?_2 $ govern the introduction of e-formulas and that $ ?\vdash_1 $ and $ ?\vdash_2$ govern their elimination. The reader may  wonder why there are two pairs of introduction and elimination rules. The reason stems from IEL's distinction between erotetic evocation and erotetic implication. Both inference types have questions as their conclusions, but while the premises of the latter also include a question, those of the former do not. To capture evocation, we need a rule that permits the introduction (elimination) of e-formulas in sequents whose antecedents (succedents) consist solely of d-wffs, while representing implication requires a rule that permits such introduction (elimination) in cases when the antecedents (succedents) also contain e-formulas. Indeed, it is natural to think that assertions license questions in a manner that differs from the way questions license further questions. In normative-pragmatic terms, the crucial difference is that entitlement to ask one question licenses the asking of another question when resolving the latter resolves the former, i.e., when the implied question promises to `move inquiry along.' Formulating two versions of each introduction and elimination rule is one way to respect this difference.
	
	Let's examine $ \vdash?_1 $. If we omit reference to background sets and side formulas, we can read the rule top-down as saying that anyone entitled to assert some statement in the set $ \{A_{1},\mydots,A_{n}\} $ is thereby entitled to ask the question $ ?\{A_{1},\mydots,A_{n}\} $, so long as she is not  entitled to assert any particular member of that set.
	
	While this rule is intended to capture erotetic evocation, it actually represents a more general relationship, since the succedents of the premise and conclusion may include a set of side formulas composed of d-wffs and/or e-formulas. The generality of the rule is motivated by the desire to formulate rules for e-formulas that conform to the standard format of those for $ \mathsf{LK} $. It also helps to satisfy cut-elimination. As we shall see, when there are no side formulas and the background set is empty, the rule faithfully depicts erotetic evocation.
	
	The first condition of the evocative relation is satisfied by the fact that the succedent of the premise, i.e.$  \{A_1,\mydots,A_n\} $,  contains all of the constitutes of the question in the succedent of the conclusion, i.e.,  $?\{A_1,\mydots,A_n\} $. Satisfying the second condition of erotetic evocation, however, requires significant exploitation of the defeat mechanism. Recall that this condition prohibits evocation whenever the set of d-wffs sc-entails a constituent of the evoked question. By including singletons of each direct answer in the defeater set of the conclusion, $ \vdash ?_1 $ guarantees that the set of auxiliary d-wffs does not sc-entail a direct answer to the evoked question. Since defeat is defined as a relation between the d-wffs in the antecedent-\textit{cum}-background set on one hand, and members of the defeater set on the other, the inclusion of direct answers as separate sets permits $ X $ to fulfill the second condition without jeopardizing its satisfaction of the first. 
	
	Example \ref{ex:simple} illustrates the dynamics of defeat for sequents involving e-formulas. In the example, $ \{p\} $ and $ \{q\} $ both belong to the defeater set of the end-sequent, thereby ensuring that neither of these formulas is derivable from the antecedent. But since $ \{p, q\} $ is not in the defeater set, the sequent is undefeated.\fnmark{empty} 
	
	\fntext{empty}{To save space, we will simply write `$ \emptyset $' in defeater sets whose composition  is technically $ \emptyset\cup\emptyset $. This notation is employed merely as a convenience; the rules remain \textit{multiplicative} not \textit{additive}.}
	
	\begin{example}\label{ex:simple}
		$\cdot\mid p\vee q\sststile{\{\{p\},\{q\}\}}{} ?\{p, q\} $
		
		$$  
		\begin{prooftree}
		\Hypo{}
		\Infer1[$ ax.$]{ \cdot\mid p\sststile{\emptyset}{}  p}
		\Hypo{}
		\Infer1[$ ax.$]{ \cdot\mid q\sststile{\emptyset}{}  q}
		\Infer2[$ \vee\vdash $]{\cdot \mid p\vee q \sststile{\emptyset}{}  p, q}
		\Infer1[$ \vdash ?_1 $]{ \cdot\mid p\vee q\sststile{\{\{p\},\{q\}\}}{}  ?\{p, q\}}
		\end{prooftree}
		$$
	\end{example}
	
	Let us now consider the right-rule for e-formulas, $ ?\vdash_1 $. Bypassing background sets and side formulas, the rule says that anyone entitled to (at least one member of) a set of assertions, $ X $, on the basis of her entitlement to assert $ A_1 $ or \ldots or $ A_n $ is also entitled to $ X $ on the basis of her entitlement to ask the question $ ?\{A_{1},\mydots,A_{n}\} $. The observant reader will no doubt have noticed that $ ?\vdash_1 $ is but an expanded version of the classical elimination rule for $ \vee $. This overlap is to be expected, since questions in IEL for classical logic \textit{presuppose} that at least one of their possible answers is true. However, $ ?\vdash_1 $ is not equivalent to $ \vee\vdash $ in $ \mathsf{LK}^\mathbf{?} $ since ,\textit{inter alia} the former is restricted to declarative succedents, while the latter is not.
	
	We now turn to the second introduction rule for e-formulas, i.e., $ \vdash ?_2 $. Ignoring background sets and side formulas, this rule may be read  top-down as follows: anyone entitled to assert $ B_1 $ or \ldots or $ B_m $ on the basis of her entitlement to ask the question $ ?\{A_{1},\mydots,A_{n}\} $ and whose entitlement to assert any member of $ \{B_{1},\mydots,B_{n}\} $ licenses her to an answer to that question, is thereby entitled to ask the question $  ?\{B_{1},\mydots,B_{n}\}  $. As noted, the rule is intended to capture the relation of regular erotetic implication. The left-most premise roughly corresponds to the first condition on R-EI, according to which an answer to the implying question, together with the set of d-wffs, entails that there is a true answer to the implied question. The premises to the right roughly correspond to the second condition, according to which, an answer to the implied question, together with the set of d-wffs, entails an answer to the implying question. Admittedly, $ \vdash ?_2 $ encodes less restrictive conditions than those imposed by R-EI insofar as it permits side formulas of either the declarative or erotetic kind to occur in the antecedent or succedent. Again, this liberalization brings the rule into conformity with standard rules in classical sequent calculi.
	
	Here is an example of a proof involving $ \vdash ?_2$.\fnmark{empty2}
	
	\fntext{empty2}{Again, we will simply write `$ \emptyset $' in defeater sets whose composition  is technically $ \emptyset\cup\emptyset $, $ \emptyset\cup\emptyset\cup \emptyset\cup\emptyset$, or $ \emptyset\cup\emptyset\cup \emptyset\cup\emptyset\cup\emptyset\cup\emptyset\cup \emptyset\cup\emptyset $.}
	
	\begin{example}
		$p, q, r\mid \: ?\{p, q\vee r\}\sststile{\emptyset}{} \: ?\{p, q, r\}$
		
		$$  
		\begin{prooftree}[separation = .5em,template = \small$\inserttext$ ]
		\Hypo{}
		\Infer1[$ \hspace{-1mm}\scriptscriptstyle ax.$]{ \cdot\mid p\sststile{\emptyset}{}  p}
		\Hypo{}
		\Infer1[$ \hspace{-1mm}\scriptscriptstyle ax.$]{\cdot\mid q\sststile{\emptyset}{}  q}
		\Hypo{}
		\Infer1[$ \hspace{-1mm}\scriptscriptstyle ax.$]{\cdot\mid r\sststile{\emptyset}{}  r}
		\Infer2[$ \hspace{-1mm}\scriptscriptstyle \vee\vdash$]{q, r\mid q\vee r\sststile{\emptyset}{}  q, r}
		\Infer2[$ \hspace{-1mm}\scriptscriptstyle ?\vdash_1$]{ q, r\mid \:?\{p, q\vee r\} \sststile{\emptyset}{}  p, q, r}
		\Hypo{}
		\Infer1[$ \hspace{-1mm}\scriptscriptstyle ax.$]{ \cdot\mid p\sststile{\emptyset}{}  p}
		\Hypo{}
		\Infer1[$ \hspace{-1mm}\scriptscriptstyle ax.$]{\cdot\mid q\sststile{\emptyset}{}  q}
		\Infer1[$ \hspace{-1mm}\scriptscriptstyle \mathsf{RW}$]{\cdot\mid q\sststile{\emptyset}{}  q, r}
		\Infer1[$ \hspace{-1mm}\scriptscriptstyle \mathsf{\vdash\vee}$]{\cdot\mid q\sststile{\emptyset}{}  q\vee r}
		\Hypo{}
		\Infer1[$ \hspace{-1mm}\scriptscriptstyle ax.$]{\cdot\mid r\sststile{\emptyset}{}  r}
		\Infer1[$ \hspace{-1mm}\scriptscriptstyle \mathsf{RW}$]{\cdot\mid r\sststile{\emptyset}{}  q, r}
		\Infer1[$ \hspace{-1mm}\scriptscriptstyle \mathsf{\vdash\vee}$]{\cdot\mid r\sststile{\emptyset}{}  q\vee r}
		\Infer4[$ \hspace{-1mm}\scriptscriptstyle \vdash?_2$]{ p, q, r\mid \: ?\{p, q\vee r\} \sststile{\emptyset}{} \: ?\{p, q, r\}}	
		\end{prooftree}
		$$
	\end{example}	
	
	Finally, we turn to $ ?\vdash_2 $. Again, ignoring background sets and side formulas, we read the rule top-down as follows: anyone who is entitled to ask the question $ ?\{B_1,\ldots,B_m\} $ on the basis of entitlement to assert any statement from the set $\{A_{1},\mydots, A_{n}\} $ and who is entitled to assert some statement in the latter on the basis of her entitlement to assert any statement in $ \{B_1,\ldots,B_m\} $, is thereby licensed to ask $ ?\{B_1,\ldots,B_m\} $ if she is entitled to ask $ ?\{A_{1},\mydots, A_{n}\} $. Recall that the motivation behind $ \vdash ?_1 $ and $ \vdash ?_2 $ is that in order to represent EE and R-EI, a calculus needs to distinguish between inferences in which a question follows from a set of statements and those in which a question follows from at least one question. Symmetry demands that we respect this difference in our elimination rules as well. Nevertheless, $ ?\vdash_2 $ and $ \vdash ?_2 $ are not on equal footing, for the former is admissible in $ \mathsf{LK^{?}} $, while the latter is not.
	
	\begin{lemma}[Admissibility of $ ?\vdash_2 $]\label{lem:?L2}
		The rule $?\vdash_2 $ is admissible in $ \mathsf{LK}^?$, i.e., if the premises of the rule are provable in $ \mathsf{LK}^?\backslash\{?\vdash_{\scriptscriptstyle 2}\}$ and its conclusion is undefeated, then the conclusion is provable in $ \mathsf{LK}^?\backslash\{?\vdash_{\scriptscriptstyle 2}\}$.
		
		\begin{proof}
			Suppose that the premises of $?\vdash_2 $ are provable.   Since the axioms of $ \mathsf{LK^{?}} $ are declarative and since the rule only stipulates that the active formula in the first set of premises be erotetic, the proof proceeds by considering the simplest case, namely, that in which the active formula $ ?[B_{\mid m}] $ is the only e-formula to occur in the premises. It follows that the first set of premises  in the rule, i.e., those of the form 
			
			\begin{equation}
			\Sigma_{i} \mid X_{i}, A_{i}\!\sststile{\mathbf{S}_{i} }{}\,?[B_{\mid m}], \varUpsilon_i
			\end{equation}
			
			\noindent for all $ i \in \{1,\ldots,n\} $, must be derived either via $ \vdash ? _1 $ from sequents of the form 
			
			\begin{equation}
			\Sigma_{i} \mid X_{i}, A_{i}\!\sststile{\mathbf{S}_{i} }{}\,[B_{\mid m}], \varUpsilon_i
			\end{equation}
			
			\noindent for all $ i \in \{1,\ldots,n\} $ (Sub-case 1) or from purely declarative sequents via $ \mathsf{RW} $ (Sub-case 2). 
			
			Sub-case 1: If (I) are obtained from (II) by $ \vdash ?_1 $, then $[\mathbb{B}_{\mid m}] \subseteq \mathbf{S}_{\mid n}$. But, since $ [B_{\mid m}] $ is added to the conclusion's background set, and since $ [B_{\mid m}] \not\succsim[\mathbb{B}_{\mid m}]$, the conclusion of $ ?\vdash_2 $ is defeated. Thus, there is no proof in the first place.
			
			Sub-case 2: If the premises of the form (I) are obtained by $ \mathsf{RW} $, from sequents of the form 
			
			\begin{equation}\label{eq:RW}
			\Sigma_{i} \mid X_i, A_i\sststile{\mathbf{S}_{i}}{} \varUpsilon_i
			\end{equation}
			
			\noindent then $ \Sigma_i \subseteq\Sigma_{\mid n+m}, X_i \subseteq \Gamma_{\mid n+m}, \mathbf{S}_i \subseteq \mathbf{S}_{\mid n}, $ and $ \Upsilon_i \subseteq \Delta_{\mid n+m}$. We thus obtain the conclusion of $ ?\vdash_2 $ by judicious applications of $ \mathsf{BE}, \mathsf{DE}, \mathsf{LW} $, and $ \mathsf{RW} $.
			
			Since the first set of premises cannot be derived via $ \vdash?_1 $, if they are not obtained by $ \mathsf{RW} $, then they must all follow by way of $ \vdash ? _2 $, and thus there is at least one erotetic side formula in the antecedent (i.e., in $ \Gamma_{\mid n}$). Suppose that $ ?[C_{\mid k}] $ is the only erotetic side formula in the antecedent of the first set of premises. It follows from our supposition that sequents of the form 
			
			\begin{equation}\label{eq:prem2?}
			\Sigma_{i} \mid X_i, ?[C_{\mid k}], A_i\sststile{\mathbf{S}_{i}}{}[B_{\mid m}], \varUpsilon_i
			\end{equation}
			
			\noindent for all $ i \in \{1,\ldots,n\} $ must be provable. We apply $ ?\vdash_1 $ to all of these sequents to obtain 
			
			\begin{equation}\label{eq:prem?2}
			\Sigma_{\mid n} \mid X_{\mid n}, ?[C_{\mid k}], ?[A_{\mid n}]\sststile{\mathbf{S}_{\mid n}}{}[B_{\mid m}], \varUpsilon_{\mid n}
			\end{equation}
			
			\noindent We then apply $ \vdash ? _2 $ to this sequent and to the remaining premises of $?\vdash_2 $ to obtain the conclusion. 
			
			In the case that there are additional erotetic side formulas in the antecedents of the first set of premises, we repeat the process of rendering sequents of the form \ref{eq:prem2?} into those of the form \ref{eq:prem?2}. Once all erotetic side formulas are thus derived, we can apply $ \vdash?_1 $ or $ \mathsf{RW} $ to obtain any erotetic side formulas in the succedent of these premises. Judicious application of $ \vdash?_1 $ and $ \vdash?_1 $ (or $ \mathsf{LW}$ and $\mathsf{RW}$) will likewise yield any erotetic side formulas in the second set of premises. We then apply $ \vdash ? _2 $ to obtain the conclusion of $?\vdash_2 $.
		\end{proof}
	\end{lemma}
	
	In the first case considered in this proof, we observe an interesting phenomena---namely, when any sequent in the first set of premises in $ ?\vdash_2 $ is derived via $\vdash ?_1$, the result is a paraproof. This occurs in virtue of the fact that the constituents of the implied question are added (as singletons) to the defeater set of $\vdash ?_1$'s conclusion and hence carry down into $ ?\vdash_2 $'s conclusion, while at the same time those constituents are added as formulas to the conclusion's background set. It is this interaction between the $ ?\vdash_2 $ and $\vdash ?_2$ that prevents the latter's admissibility. 
	
	\begin{proposition}[Non-admissibility of $\vdash ?_2$]\label{prop:?R2}
		The rule $\vdash ?_2$ is not admissible in $ \mathsf{LK^{?}} $.
		
		\begin{proof}
			To prove admissibility, we would have to derive $\vdash ?_2$'s conclusion when $ ?[A_{\mid n}] $ is the only e-formula that occurs in the rule. There are three ways to obtain the first premise (i.e., $ \Sigma_1\mid X_1,?[A_{\mid n}]\sststile{\mathbf{S}_{1}}{}[B_{\mid m}], \varUpsilon_1 $) under this condition.
			In the first case, all of the active formulas in the premise are obtained by weakening. To recover the conclusion of $ \vdash ?_2 $ we simply apply $ \mathsf{LW}, \mathsf{RW}, \mathsf{BE} $ and $ \mathsf{DE} $ as needed.
			
			In the second case, the first premise is obtained by $ \mathsf{LW} $ from a sequent of the form:
			\begin{equation}\label{eq:?}
			\Sigma_1\mid X_1 \sststile{\mathbf{S}}{} [B_{\mid m}], \varUpsilon_1.
			\end{equation}
			Since \ref{eq:?} has a declarative antecedent, we must use $ \vdash ?_1 $ to get	
			\begin{equation}\label{eq:?-}
			\Sigma_1\mid X_1 \sststile{\mathbf{S}\,\cup\,[\mathbb{B}_{\mid m}]}{} ?[B_{\mid m}], \varUpsilon_1
			\end{equation}
			Now, to arrive at the conclusion of $ \vdash ?_2 $ we must apply, \textit{inter alia} $ \mathsf{BE} $ in order to get $ [B_{\mid m}] $ into the background set. But since $ [B_{\mid m}] \not\succsim [\mathbb{B}_{\mid m}] $ the resulting sequent is defeated and there is no proof. 
			
			Finally, in the third case, the first premise is obtained by  $ ?\vdash_1 $ from a set of sequents with the following form (omitting subscripts for auxiliary sets):
			\begin{equation}\label{eq:?1}
			\Sigma\mid X, A_i \sststile{\mathbf{S}}{} [B_{\mid m}], \varUpsilon
			\end{equation}
			\noindent for all $ i \in \{1,\ldots,n\} $. We must apply $ ?\vdash_2 $ to this set to obtain the desired conclusion of $ \vdash ?_2 $. To prepare sequents of form \ref{eq:?1} to serve as premises in an application of $ ?\vdash_2 $, the e-formula $ ?[B_{\mid m}] $ must occur in their succedents. Since their antecedents are declarative, we must apply $ \vdash ?_1 $ to each of them. But as in the second case, this forces $ [\mathbb{B}_{\mid m}] $ into their defeater set. Since the conclusion of $ ?\vdash_2 $ will contain $ [B_{\mid m}]  $ in its background, the resulting sequent is defeated.
			
			Thus, unlike the proof of $ ?\vdash_2 $'s admissibility, where it was shown that all of the sequents in the first set of premises could be obtained by applying the declarative-succedent rule $ ?\vdash_1$ followed by $ \vdash ?_2 $ to yield erotetic side formulas, we cannot obtain the first premise of $ \vdash ?_2 $ by applying the declarative-antecedent rule $ \vdash ?_1$ followed by $ ?\vdash_2 $, since applying $ \vdash ?_1$ will render the desired conclusion defeated.
			
			The only instance in which the conclusion of $ \vdash ?_2 $ can be obtained from the premises in $ \mathsf{LK^{?}}\backslash \{\vdash ?_2 \} $ is when all the active formulas in the rule are obtained by weakening, i.e., the first case.
		\end{proof}
	\end{proposition}
	
	In virtue of the admissibility of $ ?\vdash_2 $ and the non-admissibility of $ \vdash?_2 $, the major results that follow will utilize $ \vdash?_2 $ whenever possible, obviating consideration of $ ?\vdash_2 $.
	
	Before moving on to our main results, we pause to prove that undefeatedness is preserved upward in \textit{all} cut-free proofs in $ \mathsf{LK}^\mathbf{?} $. As we shall see, this property plays an important role in our proofs of soundness and completeness and of cut-elimination.
	
	\begin{lemma}\label{lem:upundefeat-whole}
		Any cut-free paraproof in $ \mathsf{LK}^\mathbf{?} $ is a proof if and only if its end-sequent is undefeated, i.e., undefeatedness is preserved upwards in cut-free proofs in $ \mathsf{LK}^\mathbf{?} $.
		
		\begin{proof}
			Lemma \ref{lem:upundefeat} establishes the thesis for all cut-free proofs in $ \mathsf{LK_d}^\mathbf{?} $. The remaining proofs therefore involve application of one or more of the erotetic rules. Since there is no change to the defeater sets or antecedents in $ \vdash?_1 $, this rule obviously preserves undefeatedness upwards. As for $ ?\vdash_1 $, Lemma \ref{lem:defeatrel}.5 shows that it would be possible for the conclusion to be undefeated when one of the premises is defeated. This possibility, however, is blocked by the addition of the eliminated questions' possible answers, i.e., $ A_1, \ldots, A_n $, to the conclusion's background set, in a manner similar to $ \vee\vdash, \to\vdash, \vdash\to,\vdash\neg $. The same is done with respect to the constituents of the implied question in $ \vdash ?_2 $ and to those of both the implied and implying question in $ ?\vdash_2 $. 
		\end{proof}
	\end{lemma}
	
	\section{Soundness and Completeness Results for $ \mathsf{LK}^\mathbf{?}$}\label{sec:S+C}
	
	We will now show that a particular class of defeasible sequents in  $ \mathsf{LK}^\mathbf{?}$ is sound and compete with respect to the definitions of EE and R-EI.
	
	\begin{theorem}[Soundness with respect to EE]\label{thrm:soundness_EE}
		{\raggedright
			$\newline \text{If\:\:}   \cdot\mid X\sststile{\:[\mathbb{A}_{\mid n}] }{} ?[A_{\mid n}]  \text{\,is provable in\,\:}   \mathsf{LK}^\mathbf{?},  \text{\,then \:}   \mathbf{E}(X, ?[A_{\mid n}]).  $
		}
		
		\begin{proof}
			
			According to Definition \ref{def:erotetic_evocation}, $ \mathbf{E}(X, ?[A_{\mid n}]) $ just in case $ X\VDash [A_{\mid n}] $ and $ \neg\exists A_{i} \in [A_{\mid n}] (X \vDash A_i)  $. Since the only way to prove $ \cdot\mid X\sststile{\:[\mathbb{A}_{\mid n}] }{} ?[A_{\mid n}] $ is to apply $ \vdash ?_1 $, and since undefeatedness is preserved upwards (Lemma \ref{lem:upundefeat-whole}),  it follows by hypothesis that $ \cdot\mid X\sststile{\emptyset }{} [A_{\mid n}] $ must be provable. Hence, by Lemma \ref{lem:LKSd_sound-complete}, $ X\VDash [A_{\mid n}] $.  If $ \cdot\mid X\sststile{\:[\mathbb{A}_{\mid n}] }{} ?[A_{\mid n}] $ is provable, then it is undefeated and hence, by Lemma \ref{lem:LK_sound+complete}, $ \mathcal{E}(\emptyset\cup X)\nVDash Z $ for any $ Z\in [\mathbb{A}_{\mid n}]  $. It follows that $ X\nvDash A $ for any $ A\in[A_{\mid n}]  $. Thus, $ \neg\exists A_{i} \in [A_{\mid n}] (X \vDash A_i)  $. 
		\end{proof}
	\end{theorem}
	
	\begin{remark}
		It follows from Lemma \ref{lem:LKSd_sound} that $ \mathbf{E}(X, ?[A_{\mid n}]) $ also holds when $ \Sigma\mid X\sststile{\:[\mathbb{A}_{\mid n}] }{} ?[A_{\mid n}]$ is provable in $ \mathsf{LK^{?}} $ where $ \Sigma\neq \emptyset $.
	\end{remark}
	
	\begin{theorem}[Completeness with respect to EE]\label{thrm:completeness_EE}
		$ \newline  \text{If\:\:}  \mathbf{E}(X, ?[A_{\mid n}] ),  \text{\:then \:} \cdot\mid X\sststile{\:[\mathbb{A}_{\mid n}] }{} ?[A_{\mid n}]  \text{\,is provable in\,\:}   \mathsf{LK}^\mathbf{?}.$
		
		\begin{proof}
			If $ \mathbf{E}(X, ?[A_{\mid n}] ) $, then $ X \VDash [A_{\mid n}] $, and hence, by Lemma \ref{lem:LKSd_sound-complete}, $ \cdot\mid X\sststile{\emptyset }{} [A_{\mid n}]$ is provable. Now apply $ \vdash ?_1 $ to derive  $ \cdot\mid X\sststile{\:[\mathbb{A}_{\mid n}] }{} ?[A_{\mid n}]  $. Since the background set is empty, the result is only a paraproof if $ X\not\succsim [\mathbb{A}_{\mid n}]$. But if $ X\not\succsim [\mathbb{A}_{\mid n}]$, then $ X\VDash Z $ for some $ Z\in [\mathbb{A}_{\mid n}]$ and thus it is not the case that $ \mathbf{E}(X, ?[A_{\mid n}] ) $.
		\end{proof}
	\end{theorem}
	
	\begin{corollary}[Soundness and Completeness]\label{coro:soundness+completeness_EE}
		{\ignorespaces
			$ \newline   \mathbf{E}(X, ?[A_{\mid n}] )  \textnormal{\quad iff \quad} \cdot\mid X\sststile{\:[\mathbb{A}_{\mid n}] }{} ?[A_{\mid n}]   \text{\,is provable in\,\:}   \mathsf{LK}^\mathbf{?}. $}
		\begin{proof}
			Follows from Theorems \ref{thrm:soundness_EE} and \ref{thrm:completeness_EE}.
		\end{proof}
	\end{corollary}
	
	We now establish soundness and completeness for R-EI.
	
	\begin{theorem}[Soundness with respect to R-EI]\label{thrm:soundness_R-EI}
		
		{\ignorespaces
			$  \newline \text{If \:}  \Sigma, [ A_{\mid n}], [ B_{\mid m}] \mid X, ?[ A_{\mid n}] \sststile{\emptyset}{} \: ?[ B_{\mid m}] \text{\: is provable in \:} \mathsf{LK}^\mathbf{?}  \newline\text{for an arbitrary\,\:} \Sigma, \text{\:then\:\:} \mathbf{Im}_{r}(?[ A_{\mid n}], X, ?[ B_{\mid m}]).$}
		
		\begin{proof}
			$ \Sigma, [ A_{\mid n}], [ B_{\mid m}] \mid X, ?[ A_{\mid n}] \sststile{\emptyset}{} \: ?[ B_{\mid m}]  $ can only be obtained by weakening (Case 1) or via $ \vdash ?_2 $ (Case 2), since or $ ?\vdash_2 $ is admissible.
			
			Case 1: Consider the following proof:
			\begin{center}
				\begin{prooftree}[separation = .5em,template = \small$\inserttext$ ]
					\Hypo{}
					\Infer1[$ \scriptstyle ax.$]{ \cdot\mid p\sststile{\emptyset}{}  p}
					\Infer1[$\scriptstyle \neg \vdash $]{\cdot\mid p, \neg p\sststile{\emptyset}{}  }
					\Infer1[$\scriptstyle\mathsf{LW} $]{\cdot\mid p, \neg p, ?\{q, r\}\sststile{\emptyset}{}  }
					\Infer1[$ \scriptstyle\mathsf{RW} $]{\cdot\mid p, \neg p, ?\{q, r\}\sststile{\emptyset}{} ?\{s, t\} }
					\Infer1[$ \scriptstyle\mathsf{BE} $]{q, r, s, t\mid p, \neg p, ?\{q, r\}\sststile{\emptyset}{} ?\{s, t\} }
				\end{prooftree}
			\end{center}
			
			$ \mathbf{Im}_{r}(?\{q, r\}, \{p, \neg p\}, ?\{s, t\}) $ follows from the fact that $ A, \neg A \VDash dQ^* $ for an arbitrary $ Q^* $ and that `$ \VDash $' is monotonic.
			
			Case 2: The following must be provable:
			\begin{equation}\label{eq:alpha}
			\Sigma, [ A_{\mid n}]\mid X, ?[ A_{\mid n}] \sststile{\emptyset}{} \: [ B_{\mid m}]  
			\end{equation}
			
			Likewise, for every $ j\in \{1,\ldots, m\} $, there must exist $ A_i \in [ A_{\mid n}]  $ such that 
			
			\begin{equation}\label{eq:beta}
			\Sigma\mid X, B_j \sststile{\emptyset}{} \: A_i
			\end{equation}
			
			is provable. Thus, by Lemma \ref{lem:LKSd_sound-complete},  for each $ j\in \{1,\ldots, m\} $ there exists $ A_i\in\{A_1,\ldots,A_n\} $ such that $ X, B_j \vDash A_i $.
			
			The only way of obtaining (\ref{eq:alpha}) is by applying $ ?\vdash_1 $, since there is no other rule (save for $ \mathsf{LW} $, covered in the first case) that permits the introduction of an e-formula into an antecedent	of a sequent which has a strictly declarative succedent. So each sequent of the form:
			
			\begin{equation}
			\Sigma\mid X, A_i\sststile{\emptyset}{} \: [ B_{\mid m}]  
			\end{equation}
			
			is provable. Thus, by Lemma \ref{lem:LKSd_sound-complete} $ X, A_i \VDash \, [ B_{\mid m}] $ for each $ A_i \in [ A_{\mid n}]  $.
			
		\end{proof}
	\end{theorem}
	
	\begin{remark}
		$ \mathbf{Im}_{r}(?[ A_{\mid n}], X, ?[ B_{\mid m}]) $ also holds when $  \Sigma\mid X, ?[ A_{\mid n}] \sststile{\emptyset}{} \: ?[ B_{\mid m}] $ is provable for an arbitrary $ \Sigma $, since $  \Sigma, [ A_{\mid n}], [ B_{\mid m}] \mid X, ?[ A_{\mid n}] \sststile{\emptyset}{} \: ?[ B_{\mid m}] $ may be obtained from it via $ \mathsf{BE} $ with no threat of defeat.
	\end{remark}
	
	\begin{theorem}[Completeness with respect to R-EI]\label{thrm:completeness_R-EI}
		$ \newline \text{If \:}  \mathbf{Im}_{r}(?[ A_{\mid n}], X, ?[ B_{\mid m}]) $, then $\Sigma, [ A_{\mid n}], [ B_{\mid m}] \mid X, ?[ A_{\mid n}] \sststile{\emptyset}{} \: ?[ B_{\mid m}] \newline\text{\,is provable in }  \mathsf{LK}^\mathbf{?} \text{for an arbitrary\,\:} \Sigma.$
		
		\begin{proof}
			Suppose that $ \mathbf{Im}_{r}(?[ A_{\mid n}], X, ?[ B_{\mid m}])  $. It follows from Lemma \ref{lem:LKSd_sound-complete} that the premises of $ \vdash ?_2 $ are provable with empty defeater sets. More precisely, $ \Sigma\mid X, A_i \sststile{\emptyset}{} \: [ B_{\mid m}]  $ is provable for $ 1\leq i \leq n $, and via $ ?\vdash_1 $ we obtain $  \Sigma, [ A_{\mid n}]\mid X, ?[ A_{\mid n}]\sststile{\emptyset}{} \: [ B_{\mid m}]  $. Likewise, it follows from our supposition that for each $ j\in \{1, \ldots, m\} $ there exists $ A_i \in  [ A_{\mid n}] $ such that $ X, B_j \VDash A_i $. Hence, by Lemma \ref{lem:LKSd_sound-complete}, the corresponding sequent $ \Sigma\mid X, B_j \sststile{\emptyset}{} A_i  $ is provable. We then apply $ \vdash ?_2 $ to obtain $\Sigma, [ A_{\mid n}], [ B_{\mid m}] \mid X, ?[ A_{\mid n}] \sststile{\emptyset}{} \: ?[ B_{\mid m}] $. Since the defeater set is empty, the derivation is a proof.
		\end{proof}
	\end{theorem}
	
	\begin{corollary}[Soundness and Completeness]\label{coro:s+c_R-EI}
		$\newline \mathbf{Im}_{r}(?[ A_{\mid n}], X, ?[ B_{\mid m}])   \textnormal{\quad iff\quad} \Sigma, [ A_{\mid n}], [ B_{\mid m}] \mid X, ?[ A_{\mid n}] \sststile{\emptyset}{} \: ?[ B_{\mid m}] $.
		\begin{proof}
			Follows from Theorems \ref{thrm:soundness_R-EI} and \ref{thrm:completeness_R-EI}.
		\end{proof}
	\end{corollary}
	
	Finally, we have sequents that are sound and complete with respect to Regular Pure Erotetic Implication, i.e., RP-EI.
	
	\begin{corollary}[Soundness and Completeness]\label{coro:s+c_RP-EI}
		$\newline  \mathbf{Im}_{rpu}(?[ A_{\mid n}], ?[ B_{\mid m}])   \textnormal{\quad iff\quad} \Sigma \mid \,?[ A_{\mid n}] \sststile{\emptyset}{} \: ?[ B_{\mid m}] \text{ is provable in \,} \mathsf{LK}^\mathbf{?} \newline\text{for an arbitrary\,\:} \Sigma.$
		
		\begin{proof}
			Follows from Corollary \ref{coro:s+c_R-EI}.
		\end{proof}
	\end{corollary}
	
	\section{Cut-Elimination}\label{subsec:cut-elim}
	
	Having established the soundness and completeness of certain classes of consequence relations in $ \mathsf{LK}^\mathbf{?} $, we can now examine the global properties exhibited by the calculus, most importantly, that of cut-elimination. Naturally, the proofs of completeness (Theorems \ref{thrm:completeness_EE} and \ref{thrm:completeness_R-EI}) can be used to establish cut-elimination for the cases accounted for by these theorems. That is, they show that cut need not be used to prove sequents corresponding to erotetic evocation and implication. Our aim in this section, however, is to establish cut-elimination as a general property of $ \mathsf{LK^{?}} $.
	
	Before proceeding to our proof of cut-elimination, it is important to note that some applications of cut will never produce proofs and therefore do not require normalization.  
	
	\begin{lemma}\label{lem:cut-para}
		$ \pi $ is a paraproof if its last rule is the cut of an e-formula, Q, which is principal in applications of $ \vdash ?_1 $ and $ ?\vdash_1 $ in $ \pi$. 
		
		\begin{proof}
			Since the constituents of $ Q $ are added as a set of singletons, i.e., $ \{\{A_1\},\ldots,\{A_n\}\} $, to the defeater set of $ \vdash?_1 $'s conclusion and since the set of constituents, i.e., $ A_1,\ldots,A_n $, is added to the background of $ ?\vdash_1 $'s conclusion, an application of $ cut $ to those sequents yields the following: $ \Sigma',\Sigma, A_1,\ldots,A_n \mid \Gamma, ?\{A_1,\mydots,A_n\} \sststile{\mathbf{S}\,\cup\,\{\{A_1\},\ldots,\{A_n\}\}}{}X$. Since, e.g., $ \cdot\mid A_1,\ldots,A_n\sststile{\emptyset}{\mathcal{D}_{\mathsf{LK^{?}}}} A_1$, this sequent is defeated. Therefore, the result of any application of $ cut $ to the e-formula that is principle in the conclusions of $ \vdash?_1$ and $?\vdash_1 $ will be a paraproof.  Here is an illustration.
			
			$$  
			\begin{prooftree}[separation = 1em,template = \small$\inserttext$ ]
			\Hypo{\pi'}
			\Ellipsis{}{\Sigma_0\mid X \sststile{\mathbf{S}}{} [ A_{\mid n}]}
			\Infer1[$ \vdash ?_1$]{ \Sigma_0\mid X \sststile{\mathbf{S}\,\cup\,[ \mathbb{A}_{\mid n}]}{} ?[ A_{\mid n}]}
			\Hypo{\pi''_1}
			\Ellipsis{}{ \Sigma_1\mid \Gamma_1, A_1\sststile{\mathbf{T}_1}{}  \Upsilon\hspace{-3mm}}
			\Hypo{\pi''_n}
			\Ellipsis{}{ \ldots\:\Sigma_n\mid \Gamma_n, A_n\sststile{\mathbf{T}_n}{}  \Upsilon}
			\Infer2[$?\vdash_1$]{\Sigma_{\mid n}, [ A_{\mid n}]\mid \Gamma_{\mid n}, ?[ A_{\mid n}]\sststile{\mathbf{T}_{\mid n}}{}  \Upsilon}
			\Infer2[$ cut$]{ \Sigma_0 \cup \Sigma_{\mid n}, [ A_{\mid n}]\mid \Gamma_{\mid n}, X\sststile{\mathbf{S}\,\cup\,[ \mathbb{A}_{\mid n}]\,\cup\,\mathbf{T}_{\mid n}}{}  \Upsilon}
			\end{prooftree}
			$$
			
			Since $[ A_{\mid n}]\not\succsim  [ \mathbb{A}_{\mid n}]$, the end-sequent is defeated.
			
		\end{proof}
	\end{lemma}
	
	Once again, we observe the interaction between background and defeater set expansions that produces paraproofs in a manner previously noted in Lemma \ref{lem:?L2} and Proposition \ref{prop:?R2}.
	
	We can now proceed to our cut-elimination theorem.
	
	\begin{lemma}\label{lem:cutelim-whole}
		Any sequent that is provable in $ \mathsf{LK}^\mathbf{?} $ has a cut-free proof.
		\begin{proof}
			Lemma \ref{lem:cutelim} establishes cut-elimination for $ \mathsf{LK_d}^\mathbf{?} $. The remaining proofs therefore involve one or more of the erotetic rules. Since $ \vdash ?_1 $ and $ ?\vdash_1 $ are restricted versions of $ \vdash\vee $ and $ \vee\vdash $, respectively, we can reason from Lemma \ref{lem:cutelim} that cut-based proofs involving these rules will admit of normalization.  From Lemma \ref{lem:cut-para}, it follows that there is no need to provide a parallel reduction of $ \vdash?_1\backslash?\vdash_1 $.
			
			As $ ?\vdash_2 $ is admissible (Lemma \ref{lem:?L2}), we restrict our attention to cut-based proofs involving $ \vdash ?_2 $. There are two cases to consider. 
			
			\textbf{Case 1}: \textit{cut} is applied to a d-wff. Since d-wffs only occur as side formulas in the conclusion of $ \vdash ?_2 $, such cases are structurally analogous to cuts on side formulas of $\vee\vdash $'s conclusion and therefore their normalization follows from Lemma \ref{lem:cutelim}. 
			
			\textbf{Case 2}: The \textit{cut} rule is applied to an e-formula. There are three sub-cases depending on whether the e-formula occurs active, principle, or as a side formula in the conclusion of $ \vdash ?_2 $. 
			
			\textit{Sub-Case 2a:} Consider the derivation that results from cutting an e-formula, $ ?[ A_{\mid n}] $, that occurs principle in the succedent of $\vdash ?_1 $ and is \textit{active} in the antecedent of $ \vdash?_2 $. For the sake of concision, we introduce the following terminology: we refer to the first premise of an instance of $ \vdash ?_2 $ as its $ \alpha $ premise, i.e., $ \alpha = \Sigma_{1} \mid\Gamma_1, ?[ A_{\mid n}]\sststile{\mathbf{S}}{}[B_{\mid m}], \Delta_1 $, and refer to the remaining premises as the instance's $ \beta $ premises.   Obtaining the $ \alpha $ premise for $ \vdash?_2 $, via $ ?\vdash_1 $ requires that $ [ A_{\mid n}] $ be added to the conclusion's background set. The downward expansion of background sets entails that $ [ A_{\mid n}] $ will occur in the end-sequent's background set. From $\vdash ?_1 $, we know that $ [ \mathbb{A}_{\mid n}] $ will occur in the end-sequent's defeater set. But since $[ A_{\mid n}]\not\succsim  [ \mathbb{A}_{\mid n}]$, the end-sequent will be defeated. Again, there will be no cut-based proof. 
			
			Alternatively, the $ \alpha $ premise may be obtained by $ \mathsf{LW} $ which adds the e-formula $ ?[A_{\mid n}] $ to the antecedent without thereby bringing $ [A_{\mid n}] $ into the background set. In such a case, however, the proof must contain a sequent that does not include the cut formula. The conclusion can thus be recovered by weakening and expansion (i.e., $ \mathsf{BE} $ and \textsf{$ \mathsf{DE} $}).
			
			\textit{Sub-Case 2b:} Consider the derivations in which the principle formula of $ \vdash ?_2 $ is cut, i.e., $ ?\{B_1,\ldots,B_m\} $. There are  three sub-cases depending on whether the $ ?\{B_1,\ldots,B_m\} $ occurs active, principle, or as a side formula in the second premise of \textit{cut}. When it is active, the succedent of the second premise in \textit{cut} contains at least one e-formula, and therefore $ \vdash ?_2 $ occurs in the branch above it. We thus need a parallel reduction of $ \vdash?_2/\vdash?_2 $. To reduce the proof, we permute \textit{cut} upward and apply it to both $ \alpha $ premises and again to each of the corresponding $ \beta $ premises. Below is a simplified example of this reduction. To conserve space, we omit background sets, assume identical defeater sets for sub-proofs, and restrict the e-formulas involved to those for yes/no questions. Once again, these measures are only taken for convenience; the rules should not be construed as additive.
			$$
			\begin{prooftree}[separation = .5em,template =  {\fontsize{4.5}{4}$\inserttext$} ]
			\Hypo{\pi_1}
			\Ellipsis{}{\Gamma, ?A\sststile{\mathbf{S}}{} B,\! \neg B}
			\Hypo{\pi_2}
			\Ellipsis{}{\Gamma, B\sststile{\mathbf{S}}{} A}
			\Hypo{\pi_{3}}
			\Ellipsis{}{\Gamma, \neg B\sststile{\mathbf{S}}{} A}
			\Infer3[$\hspace{-.3em}\scriptstyle\vdash?_2$]{\Gamma, ?A\sststile{\mathbf{S}}{} ?B}
			\Hypo{\pi_5}
			\Ellipsis{}{\Gamma, B\sststile{\mathbf{T}}{} C, \neg C}
			\Hypo{\pi_6}
			\Ellipsis{}{\Gamma, \neg B\sststile{\mathbf{T}}{} C,\! \neg C}
			\Infer2[$ \hspace{-.3em}\scriptstyle ?\vdash_1 $]{\Gamma, ?B\sststile{\mathbf{T}}{} C,\! \neg C}
			\Hypo{\pi_7}
			\Ellipsis{}{\Gamma, C\sststile{\mathbf{T}}{} B}
			\Hypo{\pi_8}
			\Ellipsis{}{\Gamma, \neg C\sststile{\mathbf{T}}{} \neg B}
			\Infer3[$\scriptstyle\vdash?_2$]{\Gamma, ?B\sststile{\mathbf{T}}{} ?C}
			\Infer2[$\hspace{-.3em}\scriptstyle cut$]{\Gamma, ?A\sststile{\mathbf{S}\cup\,\mathbf{T}}{}  ?C}
			\end{prooftree}
			$$
			
			\vspace*{-4mm}
			
			$$\!\!\!\!\!\!\!\!\downsquig\!\!\!\!\!\!\!\!\!$$
			
			\vspace*{-11mm}
			
			$$
			\begin{prooftree}[separation = .5em,template =  {\fontsize{4.5}{4}$\inserttext$} ]
			\Hypo{\pi_1}
			\Ellipsis{}{\Gamma, ?A\sststile{\mathbf{S}}{} B, \!\neg B}
			\Hypo{\pi_5}
			\Ellipsis{}{\Gamma, B\sststile{\mathbf{T}}{} C,\! \neg C}
			\Infer2[$\hspace{-.5em}\scriptscriptstyle cut$]{\Gamma, ?A\sststile{\mathbf{S\cup\, T}}{} \neg B, C, \neg C}
			\Hypo{\pi_6}
			\Ellipsis{}{\Gamma, \neg B\sststile{\mathbf{T}}{} C, \neg C}
			\Infer2[$\hspace{-.5em}\scriptscriptstyle cut$]{\Gamma, ?A\sststile{\mathbf{S\cup\, T}}{} C, \neg C}
			\Hypo{\pi_7}
			\Ellipsis{}{\Gamma, C\sststile{\mathbf{T}}{} B}
			\Hypo{\pi_2}
			\Ellipsis{}{\Gamma, B\sststile{\mathbf{S}}{} A}
			\Infer2[$\hspace{-.5em}\scriptscriptstyle cut$]{\Gamma, C\sststile{\mathbf{S\cup\,T}}{} A}
			\Hypo{\pi_{3}}
			\Ellipsis{}{\Gamma, \neg B\sststile{\mathbf{S}}{} A}
			\Hypo{\pi_8}
			\Ellipsis{}{\Gamma, \neg C\sststile{\mathbf{T}}{} \neg B}
			\Infer2[$\hspace{-.5em}\scriptscriptstyle cut$]{\Gamma, \neg C\sststile{\mathbf{S\cup\,T}}{} A}
			\Infer3[$\hspace{-.5em}\scriptscriptstyle \vdash?_2$]{\Gamma, ?A\sststile{\mathbf{S}\cup\,\mathbf{T}}{}  ?C}
			\end{prooftree}
			$$

			We now consider the case when $ ?\{B_1,\ldots,B_m\} $ occurs principle in the second premise of \textit{cut}, i.e., $ \Sigma\mid \Gamma, ?\{B_1,\ldots,B_m\}\sststile{\mathbf{S}}{} \Delta  $.  This sequent must be derived via $ ?\vdash_1 $ or $ ?\vdash_2 $. But, again, as the latter is admissible, we focus on the former. Thus, sequents of the form $\Sigma_i\mid \Gamma_i, B_i\sststile{\mathbf{S}}{} \Delta $ for each $ 1\leq i \leq n $ must occur in the branch above. Since the  $ \alpha $ premise's succedent contains $ B_1,\ldots,B_m$,  reduction is achieved by iterated applications of \textit{cut}. 
			
			Lastly, if $ ?\{B_1,\ldots,B_m\} $ occurs as a side formula in the antecedent of the second premise of $ cut $, then the latter is permuted upward and applied to $ \vdash ?_2 $'s $ \alpha$ premise.
			
			\textit{Sub-Case 2c:}
			If the e-formula cut from $ \vdash ?_2 $'s conclusion is a side formula, then reduction is achieved in the fashion analogous to that of cuts on side formulas in $ \mathsf{LK} $, in keeping with Lemma \ref{lem:cutelim}.

		\end{proof}
		
	\end{lemma}

	\begin{corollary}[Subformula Property]
		If a sequent is provable in $ \mathsf{LK^?} $, then	it is provable analytically, i.e., by means of a derivation in which all formulas are subformulas of those occurring in the end sequent.
		\begin{proof}
			The result follows by induction on the length of cut-free proofs.
		\end{proof}
	\end{corollary}

	\begin{lemma}[Invertibility of the Erotetic Rules]\label{lem:erotetic_rules_invert}
		If the conclusion  of $ ?\vdash_1$,  $?\vdash_2, \vdash?_ 1, $ or $ \vdash?_ 2 $ is provable, then so is (are) its premises.
		\begin{proof}
			Straightforward by induction on proof height.
		\end{proof}
	\end{lemma}
	
	While it is not our intention  to work out the details of proof-search in $ \mathsf{LK^{?}} $, we can sketch the approach. Suppose that we wish determine whether $  \Sigma\mid\Gamma\sststile{\mathbf{S}}{}\Delta $ is provable $ \mathsf{LK^{?}} $. The first step is to determine whether the sequent is defeated, since, if it is, then it has no proof. But as defeat is defined in terms of derivability in $ \mathsf{LK^{?}} $, this step itself requires a search for derivations. To do so, we construct a set, $ \mathfrak{S} $, of sequents with the form $  \cdot\mid\mathcal{E}(\Sigma\cup\Gamma)\sststile{\emptyset}{} X_i $ for each $ X_i\in\mathbf{S} $. We then conduct a \textit{derivation search} on each of these sequents. Since the members of $ \mathfrak{S} $ are declarative sequents with empty background and defeater sets, and since the rules of $ \mathsf{LK_d}^? $ are just those of the set-based calculus $\mathsf{LK_0} $ \textit{modulo} $ \mathsf{BE}$ and $ \mathsf{DE} $, we may proceed in the manner of the standard \textit{proof search algorithm} for $ \mathsf{LK} $ (e.g. \cite[p. 295]{Bimbo2014}), with the caveat that sequents are constructed of sets and thus the rules for contraction and exchange are obviated.  If we find a derivation for a member of $ \mathfrak{S} $ then we know the original sequent is defeated and there is no proof. If the derivation search terminates without a derivation, however, then it is undefeated and, thanks to Lemma \ref{lem:upundefeat-whole}, we know that it has a proof just in case it has a derivation.  We thus conduct a derivation search on the original sequent, in keeping with the standard proof search algorithm for $ \mathsf{LK} $.\fnmark{BEDE} Further departures from the standard algorithm are attributable to the complexity of the rules for e-formulas.
	
	\fntext{BEDE}{Matters are slightly complicated by the presence of $ \mathsf{BE}$ and $ \mathsf{DE} $ in $ \mathsf{LK}^?$. However, since the background and defeater sets of any sequent are finite, we know that there is an upper bound to the number of derivations for any end-sequent containing a nonempty background or defeater set.}

	\section{The admissibility of $ \mathsf{PMC_E} $}\label{sec:PMC}
	
	As mentioned in the introduction, \cite{Wisniewski2016} has axiomatized erotetic evocation for the sequent calculus $ \mathsf{PMC_E} $. In this section, we demonstrate that the axioms and primitive rules of $ \mathsf{PMC_E} $ are respectively derivable and admissible in $ \mathsf{LK}^{\mathbf{S?}} $.
	
	As with $ \mathsf{LK}^\mathbf{?} $, $ \mathsf{PMC_E} $ is defined over the language of classical propositional logic syntactically enriched to incorporate e-formulas. But whereas $ \mathsf{LK}^\mathbf{?} $ \textit{extends} the sequent calculus for classical logic, $ \mathsf{PMC_E} $ deals with a restricted class of sequents known as e-sequents, which take the following form:
	\begin{equation}\label{eq:e-sequent}
	X \vdash ?\{A_{1},\mydots,A_{n}\}  
	\end{equation}
	where $ ?\{A_{1},\mydots,A_{n}\} $ is an e-formula as defined in $ \mathcal{L} $. As (\ref{eq:e-sequent}) illustrates, the antecedents of e-sequents are restricted to possibly empty sets of d-wffs while succedents are restricted to single e-formulas.
	
	Axioms in $ \mathsf{PMC_E} $ are defined as follows.
	
	\begin{definition}[Axioms in $ \mathsf{PMC_E} $, Complementary Literals ]\label{def:axioms_PCMe}
		
		The sequent  $ \vdash ?\{ D_{1},\mydots,D_{n} \} $ is an axiom of $ \mathsf{PMC_E} $ \textit{if and only if} (1) either each $ D_{i}  (1\leq i \leq n)$ is a literal or a disjunction of literals none of which are complementary and (2) $ D_{1}\vee\mydots\vee D_{n} $ \textit{does} involve complementary literals. Two literals are said to be \textit{complementary} if one of them is the negation of the other. 
		
	\end{definition}
	
	To prove the derivability of these axioms in $ \mathsf{LK}^\mathbf{?} $, we start with the following lemma.
	
	\begin{lemma}\label{lem:axiom_derive}
		If $ \{A_{1},\mydots,A_{n}\} $ is composed solely of complementary literals or disjunctions that do not  have complementary literals as subformulas, $ A_{1}\vee\mydots\vee A_{n}  $ does involve complementary literals, and $ n>1 $, then $ \Sigma\mid\cdot\sststile{\emptyset}{} A_{1},\mydots,A_{n} $ is provable in $ \mathsf{LK^\mathbf{?}} $, for an arbitrary $ \Sigma $.
		
		\begin{proof}
			We proceed by induction on proof height. \textbf{Base Case}: Suppose that $ \{A_{1},\mydots,A_{n}\} $ consists of two complementary literals, e.g., $ \{p, \neg p\}$. We obtain $ \cdot\mid\cdot \sststile{\emptyset}{}  p, \neg p  $ via application of $ \vdash \neg $ to the axiom introducing $ p $. 
			
			\textbf{Inductive Step}: Assume that there is a proof of $ \cdot\sststile{\emptyset}{} p, \neg p $ with height $ n $ and let $\Sigma\mid\cdot\sststile{\emptyset}{} A_{1},\mydots,A_{n} $ be the root of a proof of height $ n+1$, where $ \{A_{1},\mydots,A_{n}\} $  is a set of disjunctions that do not  involve complementary literals, but $ A_{1}\vee\mydots\vee A_{n}  $ does. It follows that at least one atomic formula is the subformula of an element of the set and its negation is that of another. Suppose that $ p $ is that atom. $ \Sigma\mid\cdot\sststile{\emptyset}{} A_{1},\mydots,A_{n} $ is derivable from $ \cdot\mid\cdot\sststile{\emptyset}{} p, \neg p $ via $ \mathsf{BE} $ and the rules of $ \mathsf{LK_d}^? $. 
			
		\end{proof}
	\end{lemma}
	
	\begin{theorem}\label{thrm:PCMe_axiom_derive}
		If $ \vdash ?\{A_{1},\mydots,A_{n}\} $ is an axiom in $ \mathsf{PMC_E} $, then \newline$ \Sigma\mid\cdot\sststile{[\mathbb{A}_{\mid n}]}{} ?\{A_{1},\mydots,A_{n}\} $ is provable in $ \mathsf{LK^\mathbf{?}} $ when $ \Sigma\succsim [\mathbb{A}_{\mid n}] $.
		\begin{proof}
			Lemma \ref{lem:axiom_derive} establishes that $ \Sigma\mid\cdot\sststile{\emptyset}{} A_{1},\mydots,A_{n} $ is provable in $ \mathsf{LK^\mathbf{?}} $, where $ \{A_{1},\mydots,A_{n}\} $ is composed solely of complementary literals or disjunctions that do not  have complementary literals as subformulas, $ A_{1}\vee\mydots\vee A_{n}  $ does involve complementary literals, and $ n>1 $. $ \Sigma\mid\cdot\sststile{[\mathbb{A}_{\mid n}]}{} ?\{A_{1},\mydots,A_{n} \}$ is derivable from $ \Sigma\mid\cdot\sststile{\emptyset}{} A_{1},\mydots,A_{n} $ via $ \vdash ?_1 $. Since by hypothesis, $ \Sigma\succsim [\mathbb{A}_{\mid n}] $, $ \Sigma\mid\cdot\sststile{[\mathbb{A}_{\mid n}]}{} ?\{A_{1},\mydots,A_{n} \}$ is undefeated and therefore provable. The thesis then follows according to Definition \ref{def:axioms_PCMe}.
		\end{proof}
	\end{theorem}
	
	\begin{definition}[Primitive Rules of $ \mathsf{PMC_E} $]\label{{def:PMCe_prim_rules}}
		
		There are four primitive rules in $ \mathsf{PMC_E} $. Note that $ A_i $ denotes any arbitrary formula in  $ \{A_1,\mydots,A_n\} $.
		\vspace{.5cm}
		\begin{itemize}
			\setlength\itemsep{1em}
			
			\item[$ R_1 $]	
			\begin{prooftree}
				\Hypo{X\vdash ?\{A_1,\mydots,A_n, B\}}
				\Hypo{X\vdash ?\{A_1,\mydots,A_n, C\}}
				\Infer2[\footnotesize$\:\text{Provided that}\: (B\wedge C) \neq A_i $.]{ X\vdash ?\{A_1,\mydots,A_n, B\wedge C\}}
			\end{prooftree}
			
			\item[$ R_2 $]
			\begin{prooftree}
				\Hypo{X\vdash ?\{A_1,\mydots,A_n, B\}}
				\Infer1[\footnotesize$\text{where}\:(B\leftrightarrow C)\:\text{ is a theorem of CPL and}\:  C \neq A_i $.]{ X\vdash ?\{A_1,\mydots,A_n,  C\}}
			\end{prooftree}
			
			\item[$ R_3 $]
			\begin{prooftree}
				\Hypo{X\vdash ?\{B \to A_1,\mydots,B\to A_n\}}
				\Infer1[]{ X, B \vdash ?\{A_1,\mydots,A_n\}}
			\end{prooftree}
			
			\item[$ R_4 $]
			\begin{prooftree}
				\Hypo{X\vdash ?\{A_1,\mydots,A_n\}}
				\Infer1[\footnotesize$\text{where}\:\mathbf{d}?\{A_1,\ldots,A_n\} \,= \,\mathbf{d}?\{B_1,\ldots,B_m\}$.]{ X\vdash ?\{B_1,\mydots,B_m\}}
			\end{prooftree}
		\end{itemize}
	\end{definition}
	
	\vspace{2mm}
	
	We proceed to demonstrate that these rules are admissible in $ \mathsf{LK^{?}} $. Since sequents in $ \mathsf{LK^{?}} $ contain background and defeater sets which are not present in the sequents of $ \mathsf{PMC_E} $, the admissibility of these rules needs to be clarified. 
	
	\begin{definition}[$ \mathsf{LK^{?}} $-correlate]
		Let  $ \mathcal{S} $ be a sequent, $\Sigma\mid X \sststile{\mathbf{S}}{} \, ?\{A_1,\mydots,A_n\} $, such that if it is provable in $ \mathsf{LK^{?}} $ then there is \textit{no} sequent $\mathcal {S'}=\Sigma'\mid X \sststile{\mathbf{S'}}{} \, ?\{A_1,\mydots,A_n\} $ provable in $ \mathsf{LK^{?}} $ for which either $ \Sigma'\subset\Sigma $ or $ \mathbf{S'}\subset\mathbf{S} $. $ \mathcal{S} $ is said to be the $ \mathsf{LK^{?}} $-correlate of the sequent $ X\vdash ?\{A_1,\mydots,A_n,\} $ in $ \mathsf{PMC_E} $. 
	\end{definition}
	
	While the background and defeater sets of $ \mathsf{LK^{?}} $-correlates are not specified, neither are they arbitrary. Rather, they should be the smallest sets with which a sequent composed of its antecedent and succedent is provable in $ \mathsf{LK^{?}} $. We may safely assume that proofs for $ \mathsf{LK^{?}} $-correlates, if they exist, do not make use of $ \mathsf{BE} $ or $ \mathsf{DE} $.
	
	\begin{definition}[Admissibility of $ \mathsf{PMC_E} $ rules in $ \mathsf{LK^{?}} $]
		A rule $ \mathcal{R} $ of $ \mathsf{PMC_E} $ is admissible in $ \mathsf{LK^{?}} $ iff there is a provable $ \mathsf{LK^{?}} $-correlate of its conclusion if there is a provable $ \mathsf{LK^{?}} $-correlate of its premise(s).
	\end{definition}

	\begin{lemma}[Admissibility of $ R_1 $]\label{lem:R1_admissible}
		The primitive rule $ R_1 $ of $ \mathsf{PMC_E} $ is admissible in $ \mathsf{LK^\mathbf{?}} $.
		
		\begin{proof}
			Assume that $ \Sigma'\mid X \sststile{[\mathbb{A}_{\mid n}]\cup\{\{B\}\} }{} \, ?\{A_1,\mydots,A_n, B\} $ and $ \Sigma''\mid X \sststile{[\mathbb{A}_{\mid n}]\cup\{\{C\}\} }{} \, ?\{A_1,\mydots,A_n, C\} $ are provable $ \mathsf{LK^{?}}$-correlates of the premises in $ R_1 $. From Lemma \ref{lem:erotetic_rules_invert}, it follows that $ \Sigma'\mid X \sststile{\emptyset }{} A_1,\mydots,A_n, B $ and $ \Sigma''\mid X \sststile{\emptyset}{} A_1,\mydots,A_n, C $ are provable. $ \Sigma\mid X \sststile{\emptyset}{} A_1,\mydots,A_n, B\wedge C  $ follows via $ \vdash\wedge $. By applying $ \vdash ?_1 $, we obtain $ \Sigma\mid X \sststile{[\mathbb{A}_{\mid n}]\cup\{\{B\wedge C\}\} }{} ?\{A_1,\mydots,A_n, B\wedge C\} $, the $ \mathsf{LK^{?}}$-correlate of the conclusion in $ R_1$ .   
		\end{proof}
	\end{lemma}
	
	\begin{lemma}[Admissibility of $ R_2 $]\label{lem:R2_admissible}
		The primitive rule $ R_2 $ of $ \mathsf{PMC_E} $ is admissible in $ \mathsf{LK^\mathbf{?}} $.
		
		\begin{proof}
			Assume that $ \Sigma\mid X\sststile{[\mathbb{A}_{\mid n}]\cup\{\{B\}\} }{} ?\{A_1,\mydots,A_n, B\} $ is a provable $ \mathsf{LK^{?}}$-correlate of the premise in $ R_2 $. From Lemma \ref{lem:erotetic_rules_invert} it follows that $ \Sigma\mid X\sststile{\emptyset}{} A_1,\mydots,A_n, B$ is provable.  Let $B\leftrightarrow C$ be a theorem of CPL where $C \neq A_i $. From Lemma \ref{lem:LKSd_sound-complete}, it follows that $ \Sigma'\mid X\sststile{\emptyset}{} A_1,\mydots,A_n, C$ is provable (where it may be that $ \Sigma'\neq\Sigma $). We apply $ \vdash ?_1 $ to obtain $ \Sigma'\mid X\sststile{[\mathbb{A}_{\mid n}]\cup\{\{C\}\} }{} A_1,\mydots,A_n, C$, the $ \mathsf{LK^{?}}$-correlate of the conclusion in $ R_2 $.
		\end{proof}
	\end{lemma}
	
	\begin{lemma}[Admissibility of $ R_3 $]\label{lem:R3_admissible}
		The primitive rule $ R_3 $ of $ PMC_E $ is admissible in $ \mathsf{LK^\mathbf{?}} $.
		
		\begin{proof}
			Suppose that $ \Sigma, B\mid X \sststile{\mathbf{S} }{} ?\{B\to A_1,\mydots,B\to A_n\} $ is provable in $ \mathsf{LK^\mathbf{?}} $ where $ \mathbf{S} = \{\{B\to A_1\},\ldots,\{B\to A_n\}\} $. From Lemma \ref{lem:erotetic_rules_invert}, it follows that $ \Sigma, B\mid X \sststile{\emptyset}{} B\to A_1,\mydots,B\to A_n$ is provable. From $ \vdash \to $ it follows that $ \Sigma\mid X, B \sststile{\emptyset}{}  A_i, \ldots, A_j$ for $ 1\leq i, j \leq n $ (where it may be that $ i=j $). If needed, we apply $ \mathsf{RW} $ to obtain $ \Sigma\mid X, B \sststile{\emptyset}{}  A_1,\ldots,A_n$ and apply $ \vdash ?_1 $ to the latter to obtain $ \Sigma\mid X, B \sststile{\mathbf{S} }{}  ?\{A_1,\mydots, A_n\}$, the $ \mathsf{LK^{?}}$-correlate of the conclusion in $ R_3 $.
		\end{proof}
	\end{lemma}
	
	\begin{theorem}[Admissibility of Primitive Rules of $ PMC_E$]\label{thrm:R-rules_admissible}
		The primitive rules $ R_1, R_2, R_3, R_4 $ of $ PMC_E $ are admissible in $ \mathsf{LK^\mathbf{?}} $.
		
		\begin{proof}
			The admissibility of $ R_1, R_2 $ and $ R_3 $ is established by Lemmas \ref{lem:R1_admissible}, \ref{lem:R2_admissible} and \ref{lem:R3_admissible}. The admissibility of $ R_4 $ is straightforward and its proof is left as an exercise for the reader. 
		\end{proof}
	\end{theorem}
	
	To our eyes, $ \mathsf{LK^{?}} $ has several advantages over $ \mathsf{PMC_E} $. The most obvious is that $ \mathsf{LK^{?}} $ encodes rules for erotetic implication and that its rules for erotetic evocation capture that relation's defeasible character. But we also believe that it is considerably easier to construct proofs in $ \mathsf{LK^{?}} $, e.g., Example \ref{ex:simple}, than to conduct proof-search in $ \mathsf{PMC_E} $, since the latter deploys axioms that must meet complicated syntactic constraints. Moreover, the proviso on $ R_2 $ shows that  $ \mathsf{PMC_E} $ must appeal to external considerations---i.e., the theorems of CPL---in order to construct proofs, whereas $ \mathsf{LK^{?}} $ is a self-contained and self-sufficient calculus. Finally, $ \mathsf{LK^{?}} $ permits sequents with multiple e-formulas occurring on either side of the turnstile, while sequents in $ \mathsf{PMC_E} $ can only contain one e-formula and only in the succedent. For this reason, $ \mathsf{LK^{?}} $ appears to have greater expressive power.

	\section{Conclusion}\label{sec:conclusion}
	
	This paper has sought to advance the formal treatment of erotetic inferences by providing a proof theory for defeasible and transitive species of the inferences studied by IEL. We have constructed a sequent calculus that (1) encodes defeasibility at the level of the turnstile; that (2) contains classes of consequences relations which are sound and complete for \textit{erotetic evocation} (Corollary \ref{coro:soundness+completeness_EE}), \textit{regular erotetic implication} (Corollary \ref{coro:s+c_R-EI}), and \textit{regular pure erotetic implication} (Corollary \ref{coro:s+c_RP-EI}); and that (3) provides introduction and elimination rules for erotetic formulas. When combined with the fact that the rules of $ PMC_E$ are admissible in our calculus, these results demonstrate that our proof theory surmounts the limitations of \cite{Wisniewski2016}'s system. 
	
	One way in which $ \mathsf{LK^\mathbf{?}} $ might be fruitfully developed would be to extend it in a such a way as to capture the general relation of erotetic implication, i.e., $ \mathbf{Im}(Q, X, Q_1) $ (Definition \ref{def:e-imp}). To do so, we might let $ \Theta_i, \Theta_j $ stand for nonempty proper subsets of $ [ A_{\mid n}] $, i.e., $ \Theta_i, \Theta_j \in \mathcal{P}([ A_{\mid n}])/\{\emptyset, [ A_{\mid n}]\} $. By replacing $ A_i, A_j $ in $ \vdash?_2 $ with $ \Theta_i, \Theta_j   $, respectively, it may be possible to obtain a rule that yields sequents that are sound and complete with respect to erotetic implication. Similar changes would need to be made to $ ?\vdash_2 $.
	
	Unfortunately, such an extension could not represent the relation of \textit{strong erotetic implication} ($ \mathbf{Im}_{s}(Q, X, Q_1)  $), a restriction of erotetic implication to those instances in which the set of declaratives does not mc-entail a nonempty proper subset of direct answers to the implying question. Strong erotetic implication is thus defeasible in a manner similar to erotetic evocation. To formulate rule(s) corresponding to $ \mathbf{Im}_{s}(?[A_{\mid n}], X, ?[ B_{\mid m}]) $ we would need to add all nonempty proper subsets of constituents of the implying question to the conclusion(s)'s defeater set---i.e.,  $\mathcal{P}([ A_{\mid n}])/\{\emptyset, [ A_{\mid n}]\}\subseteq \mathbf{S} $. Note, however, that the conclusion of $?\vdash_2$ contains $ [ A_{\mid n}] $ in its background set and that the first premise of $ \vdash ?_2 $, whether obtained by  $ ?\vdash_1 $ or $?\vdash_2$, will also have $ [ A_{\mid n}] $ in its background set. Since $ [ A_{\mid n}] \not\succsim\mathcal{P}([ A_{\mid n}])/\{\emptyset, [ A_{\mid n}]\}$, any such rule for strong erotetic implication would uniformly yield praraproofs. So capturing strong erotetic implication would require substantial modification to the calculus.
	
	Since we have sketched the contours of proof search in $ \mathsf{LK^{?}} $, another direction of development would be to construct a tool for finding proofs. Given the computational benefits of a cut-free sequent system, the development of such a tool would not be terribly difficult. We invite future research in this area.
	
	Finally, we take our results as indicating the fruitfulness of normative-pragmatic interpretations of erotetic inferences. Since the normative-pragmatic rendering for inferences is associated with proof-theoretic or inferentialist approaches to semantics, we would be interested to see whether $ \mathsf{LK^\mathbf{?}} $, or some similar system, might serve as the basis for a proof-theoretic semantics of questions. As such approaches rely on the formulation of introduction and/or elimination rules for logical connectives, the fact that our calculus provides these rules for question-forming expressions suggests that such an endeavor has some \textit{prima facie} plausibility. Once again, we invite future research in this area.

	\renewcommand{\thefigure}{1}
	\begin{landscape}
		\begin{figure}[!ht]  
			\caption{Rules for $\mathsf{LK}^\mathbf{?} $}

			\hspace*{-2cm}
			\textbf{Logical Axioms}
			\hspace{4.5cm}
			\textbf{Cut Rule}
			\vspace{.3cm}
			
			\hspace{2mm}
			$ 	\begin{prooftree}
			\Hypo{}
			\Infer1[$ax.$]{ \cdot\mid p\sststile{\emptyset}{} p}
			\end{prooftree} 
			\hspace{4.2cm}
			\begin{prooftree}
			\Hypo{\Sigma\mid\Gamma\sststile{\mathbf{S}}{}F,\Delta}
			\Hypo{\Sigma'\mid\Gamma', F\sststile{\mathbf{T}}{}\Delta'}
			\Infer2[$cut$]{ \Sigma',\Sigma\mid\Gamma',\Gamma\sststile{\mathbf{S}\,\cup\,\mathbf{T}}{}\Delta,\Delta'}
			\end{prooftree} $
			
			\vspace{.3cm}
			\textbf{Structural Rules}
			\vspace{.3cm}
			
			\hspace{2mm}
			$ 	\begin{prooftree}
			\Hypo{\Sigma\mid\Gamma\sststile{\mathbf{S}}{}\Delta}
			\Infer1[ $ \mathsf{LW} $ ]{ \Sigma\mid\Gamma, F\sststile{\mathbf{S}}{}\Delta}
			\end{prooftree}
			\hspace{2cm}
			\begin{prooftree}
			\Hypo{\Sigma\mid\Gamma\sststile{\mathbf{S}}{}\Delta}
			\Infer1[ $ \mathsf{RW} $ ]{ \Sigma\mid\Gamma\sststile{\mathbf{S}}{}F, \Delta}
			\end{prooftree} 
			\hspace{2cm}
			\begin{prooftree}
			\Hypo{\Sigma\mid\Gamma\sststile{\mathbf{S}}{}\Delta}
			\Infer1[ $ \mathsf{DE} $ ]{ \Sigma\mid\Gamma\sststile{\mathbf{S}\,\cup\,\mathbf{T}}{}\Delta}
			\end{prooftree}
			\hspace{2cm}
			\begin{prooftree}
			\Hypo{\Sigma\mid\Gamma\sststile{\mathbf{S}}{}\Delta}
			\Infer1[ $ \mathsf{BE} $ ]{ \Sigma, F\mid\Gamma\sststile{\mathbf{S}}{}\Delta}
			\end{prooftree} $
			
			\vspace{.3cm}
			\textbf{Logical Rules for D-wffs}
			
			\hspace{2mm}
			$ 		\begin{prooftree}
			\Hypo{\Sigma\mid\Gamma,A,B\sststile{\mathbf{S}}{}\Delta}
			\Infer1[$\wedge\vdash$]{ \Sigma\mid\Gamma,A\wedge B\sststile{\mathbf{S}}{}\Delta}
			\end{prooftree}
			\hspace{.5cm}
			\begin{prooftree}
			\Hypo{\Sigma\mid\Gamma\sststile{\mathbf{S}}{}A,\Delta}
			\Hypo{\Sigma'\mid\Gamma'\sststile{\mathbf{T}}{}B,\Delta'}
			\Infer2[$\vdash\wedge$]{ \Sigma',\Sigma\mid\Gamma',\Gamma\sststile{\mathbf{S}\,\cup\, \mathbf{T}}{}A\wedge B,\Delta,\Delta'}
			\end{prooftree} 
			\hspace{.5cm}		
			\begin{prooftree}
			\Hypo{\Sigma\mid\Gamma,A\sststile{\mathbf{S}}{}\Delta}
			\Hypo{\Sigma'\mid\Gamma',B\sststile{\mathbf{T}}{}\Delta'}
			\Infer2[$ \vee\vdash$]{ \Sigma',\Sigma, A, B\mid\Gamma',\Gamma, A\vee B\sststile{\mathbf{S}\,\cup\,\mathbf{T}}{}\Delta,\Delta'}
			\end{prooftree}
			\hspace{.5cm}
			\begin{prooftree}
			\Hypo{\Sigma\mid\Gamma\sststile{\mathbf{S}}{}A,B,\Delta}
			\Infer1[$\vdash\vee$]{ \Sigma\mid\Gamma\sststile{\mathbf{S}}{}A\vee B,\Delta}
			\end{prooftree}$
			
			\vspace{.3cm}
			\hspace{2mm}
			$
			\begin{prooftree}
			\Hypo{\Sigma\mid\Gamma\sststile{\mathbf{S}}{}A,\Delta}
			\Hypo{\Sigma'\mid\Gamma',B\sststile{\mathbf{T}}{}\Delta'}
			\Infer2[$\to\vdash$]{ \Sigma',\Sigma,\neg A,B\mid\Gamma',\Gamma,A\to B\sststile{\mathbf{S\,\cup\, T}}{}\Delta,\Delta'}
			\end{prooftree}
			\hspace{.5cm}	
			\begin{prooftree}
			\Hypo{\Sigma\mid\Gamma,A\sststile{\mathbf{S}}{}B,\Delta}
			\Infer1[$\vdash\to$]{ \Sigma,A\mid\Gamma\sststile{\mathbf{S}}{}A\to B,\Delta}
			\end{prooftree}
			\hspace{.5cm}
			\begin{prooftree}
			\Hypo{\Sigma\mid\Gamma\sststile{\mathbf{S}}{}A,\Delta}
			\Infer1[$\neg\vdash$]{ \Sigma\mid\Gamma,\neg A\sststile{\mathbf{S}}{}\Delta}
			\end{prooftree}
			\hspace{.5cm}
			\begin{prooftree}
			\Hypo{\Sigma\mid\Gamma,A\sststile{\mathbf{S}}{}\Delta}
			\Infer1[$\vdash\neg$]{ \Sigma,A\mid\Gamma\sststile{\mathbf{S}}{}\neg A, \Delta}
			\end{prooftree} $
			
			\vspace{.3cm}

			\textbf{Logical Rules for E-formulas}
			
			
			\begin{center}
				$\hspace{2mm}
				\begin{prooftree}
				\Hypo{\Sigma\mid X\sststile{\mathbf{S} }{} A_1,\ldots,A_n, \Delta}
				\Infer1[$\vdash?_{\scriptscriptstyle 1}\,^\dag$]{\Sigma\mid X\sststile{\mathbf{S}\,\cup\,\{\{A_1\},\ldots,\{A_n\}\}}{}
					\: ?\{A_1,\ldots,A_n\}, \Delta}
				\end{prooftree}
				\hspace{.5cm}
				\begin{prooftree}
				\Hypo{\Sigma_1\mid \Gamma_1, A_1 \sststile{\mathbf{S}_1}{} X}
				\Hypo{\hspace{-5.5mm}\ldots\Sigma_n\mid \Gamma_n, A_n \sststile{\mathbf{S}_n }{} X}
				\Infer2[$?\vdash_{\scriptstyle 1}\,\!\!\!^\dag$]{\Sigma_{\mid n}, A_1,\ldots,A_n \mid \Gamma_{\mid n}, ?\{A_1,\ldots,A_n\} \sststile{\mathbf{S}_1\cup\ldots\cup\,\mathbf{S}_{n}}{}X}
				\end{prooftree}$
				
				\vspace{.4cm}
				
				\begin{prooftree}
					\Hypo{\Sigma_{1} \mid\Gamma_1, ?[A_{\mid n}]\sststile{\mathbf{S}}{}[B_{\mid m}], \Delta_1}
					
					\Hypo{\Sigma_{2} \mid\Gamma_2, B_{1}\sststile{\mathbf{T}_{1} }{} A_i, \Delta_2}
					\Hypo{\hspace{-5mm}\ldots\Sigma_{m+1} \mid\Gamma_{m+1}, B_{m}\sststile{\mathbf{T}_{m} }{} A_{j}, \Delta_{m+1}}
					\Infer3[$\vdash?_{\scriptstyle 2}\,^\ddag$]{ \Sigma_{\mid m+1}, [B_{\mid m}]\mid\Gamma_{\mid m+1} , ?[A_{\mid n}] \sststile{\mathbf{S}\,\cup\,\mathbf{T}_{\mid m}}{}  \,?[B_{\mid m}], \Delta_{\mid m+1}}
				\end{prooftree}
				
				\vspace{.4cm}
				
				\begin{prooftree}
					\Hypo{\Sigma_{1} \!\!\mid\!\Gamma_{1}, A_{1}\!\sststile{\mathbf{S}_{1} }{}\,?[B_{\mid m}], \Delta_1 \!\ldots}
					\Hypo{\hspace{-5mm}\Sigma_{n}\!\! \mid\!\Gamma_{n}, A_{n}\!\sststile{\mathbf{S}_{n} }{}\,?[B_{\mid m}], \Delta_n }
					\Hypo{\Sigma_{n+1} \!\!\mid\!\Gamma_{n+1}, B_{1}\!\sststile{\mathbf{T}_{1} }{} A_i , \Delta_{n+1} \!\ldots}
					\Hypo{\hspace{-5mm}\Sigma_{n+m}\!\! \mid\!\Gamma_{n+m}, B_{m}\!\sststile{\mathbf{T}_{m} }{} A_j, \Delta_{n+m}}
					\Infer4[$?\vdash_{\scriptstyle 2}\,^{\!\!\!\ddag}$]{\Sigma_{\mid n+m}, [A_{\mid n}], [B_{\mid m}]\mid\Gamma_{\mid n+m}, ?[A_{\mid n}] \sststile{\mathbf{S}_{\mid n}\cup\,\mathbf{T}_{\mid m}}{}\, ?[B_{\mid m}], \Delta_{\mid n+m}}
				\end{prooftree}
				
			\end{center}
			
			\vspace*{3mm}

				$ \dag $ Provided $ n>1 $ and $ A_1,\ldots,A_n $ are non-equiform.
				
				$ \ddag $ Provided $A_{i},A_{j} \in \{ A_1,\ldots,A_n \}$; $ m,n>1$ ;  $ A_1,\ldots,A_n $ are pairwise non-equiform; and  $ B_1,\ldots,B_n $ are pairwise non-equiform.
			
		\end{figure}
	\end{landscape}

	%
	%
	%
\end{document}